\documentclass[12pt]{article}
\usepackage{amssymb,amsfonts}
\usepackage{amscd}
\usepackage{verbatim}

\newtheorem{theorem}{Theorem}[section]
\newtheorem{lemma}[theorem]{Lemma}
\newtheorem{proposition}[theorem]{Proposition}

\newtheorem{example}[theorem]{Example}
\newtheorem{definition}[theorem]{Definition}
\newtheorem{corollary}[theorem]{Corollary}
\newtheorem{remark}[theorem]{Remark}

\newenvironment{proof}{\bf Proof. \rm}{$\Box$}

\newcommand{\be}{\begin{equation}}
\newcommand{\ee}{\end{equation}}

\newcommand{\norm}[1]{\Vert #1 \Vert}

\begin{document}
\title{Representations of product systems over semigroups and dilations of commuting CP maps}

\author{ Baruch Solel\thanks{Supported by the Fund for the Promotion of Research at the Technion.}
\\Department of Mathematics\\Technion\\32000 Haifa, Israel
\\e-mail: mabaruch@techunix.technion.ac.il}

\date{}

\maketitle

\begin{abstract}
 We
prove that every pair of commuting CP maps on a von Neumann
algebra $M$ can be dilated to a commuting pair of endomorphisms
(on a larger von Neumann algebra). To achieve this,
 we first prove that every completely contractive  representation of
 a product system of $C^*$-correspondences over the
semigroup $\mathbb{N}^2$ can be dilated to an isometric (or
Toeplitz) representation.
\textbf{ key words and phrases.} completely positive maps,
 correspondences, product systems, endomorphisms, dilations.

\end{abstract}

\maketitle
\begin{section}{Introduction}
A $C^*$-correspondence $E$ over a $C^*$-algebra $A$ is a (right)
Hilbert $C^*$-module over $A$ that carries also a left action of
$A$ (by adjointable operators). It is also called a Hilbert
bimodule in the literature. A c.c. representation of $E$ on a
Hilbert space $H$ is a pair $(\sigma,T)$ where $\sigma$ is a
representation of $A$ on $H$ and $T:E\rightarrow B(H)$ is a
completely contractive linear map that is also a bimodule map
(that is, $T(a\cdot \xi \cdot b)=\sigma(a)T(\xi)\sigma(b)$ for
$a,b \in A$ and $\xi \in E$). The representation is said to be
isometric (or Toeplitz ) if $T(\xi)^*T(\eta)=\sigma(\langle
\xi,\eta \rangle)$ for every $\xi,\eta \in E$.

 In \cite{P}, Pimsner associated with such a correspondence two
 $C^*$-algebras ($\mathcal{O}(E)$ and $\mathcal{T}(E)$) with certain universal
  properties. In \cite{MS98} we studied a universal operator
  algebra (called the tensor algebra) $\mathcal{T}_+(E)$
  associated with such a correspondence.

A product system $X$ of $C^*$-correspondences over a semigroup $P$
is, roughly speaking, a family $\{X_s: s\in P\}$ of
$C^*$-correspondences
(over the same $C^*$-algebra $A$), with $X_e=A$, such that $X_s \otimes X_t$ is
isomorphic to $X_{st}$ for all $s,t \in P\backslash \{e\}$  (See Section~\ref{ps}
for the precise definition). A c.c. (respectively, isometric)
 representation of $X$ is a family
$\{T_s\}$ such that, for all $s\in P$, $(T_e,T_s)$ is
a c.c. (respectively, isometric) representation of $X_s$ for all $s\in
P$ and such that, whenever $x\in X_s$ and $y\in X_t$,
$T_{st}(\theta_{s,t}(x\otimes y))=T_s(x)T_t(y)$ (where
$\theta_{s,t}:X_s \otimes X_t \rightarrow X_{st}$ is the
isomorphism).

 If $A=\mathbb{C}$, a $C^*$-correspondence over $A$ is simply a
 Hilbert space. In \cite{Arv1}, Arveson introduced product systems of
 Hilbert spaces over the semigroup $\mathbb{R}_+$ (in order to study semigroups of
of endomorphisms of $B(H)$). When the semigroup is not discrete,
one usually assumes certain continuity or measurability conditions
on the product system.
Product systems of $C^*$-correspondences over $\mathbb{R}_+$ or
subsemigroups of $\mathbb{R}_+$ were studied by various authors
(e.g. \cite{BS},\cite{Di},\cite{MSQMP},\cite{Ma},\cite{Sk},\cite{H} and
others). Of course, a single correspondence can also be thought of
as a product system over the semigroup $\mathbb{N}$.

In \cite{F}, Fowler studied product systems over more general
(discrete) semigroups $P$. He proved the existence of a
$C^*$-algebra $\mathcal{T}(X)$ that is universal with respect to
Toeplitz representations. In fact, in most of the work done on
operator algebras associated with product systems of
correspondences (on semigroups other than $\mathbb{N}$),
 the operator algebras that were studied are
$C^*$-algebras. Two exceptions that we are aware of are \cite{KP}
and \cite{Dpo}.
 The fast growing body of litrature dealing with
$C^*$-algebras associated with $k$-graphs (see \cite{R} and the
references there) can also be viewed as the study of certain
product systems of correspondences over the semigroup
$\mathbb{N}^k$. Among other works on product systems over
semigroups, see \cite{F99},\cite{FR},\cite{MSY} and \cite{La}.

In Section~\ref{ps} we associate, with every product system $X$ of
$C^*$-\\correspondences over a discrete semigroup $P$ (with a unit
and left-cancellation), an operator algebra $\mathcal{T}_+(X)$
(called the \emph{universal tensor algebra}) which is universal
with respect to completely contractive representations. (See
Proposition~\ref{universal}).

For the rest of the paper (Sections 4 and 5) we concentrate on the
case where the semigroup $P$ is $\mathbb{N}^2$.

One of the main results of the paper (Theorem~\ref{enddilation})
 is a dilation result for a
pair of commuting (contractive, normal) completely positive  maps
on a von Neumann algebra $M$ (to a pair of commuting normal
$^*$-endomorphisms on a larger von Neumann algebra $R$). A special
case ($M=B(H)$) was proved by Bhat in \cite{Bh2} but the methods
used here are very different and the emphasis here is on the
relationship between representations of product systems and
semigroups of CP maps (as explained below).

Over the years there have been numerous studies wherein the
authors dilate CP maps or semigroups of CP maps. One can find in
the literature several approaches to dilation theory (of semigroups of CP maps)
 with
different properties. For a recent account and a list of
references see \cite[Chapter 8]{Arv2}. We shall concentrate here on
the dilations of the kind that proved effective in the study of
CP-semigroups and $E_0$-semigroups initiated by Powers and
Arveson.

Suppose $M$ is a von Neumann algebra acting on a Hilbert space
 $H$ and $\Theta$ is a contractive, normal, completely positive
 map on $M$. A quadruple $(K,R,\alpha,W)$, consisting of a Hilbert
 space $K$, a von Neumann algebra $R$, a $^*$-endomorphism
 $\alpha$ and an isometric embedding $W$ of $H$ into $K$ will be
 called an \emph{endomorphic dilation} of $(M,\Theta)$ if
 $\alpha(WW^*)WW^*=\alpha(I)WW^*$ (i.e., $WW^*$ is coinvariant
 under $\alpha$), $W^*RW=M$ (i.e. $M$ embeds as a corner of $R$)
 and, for all $a\in M$,
 $$\Theta(a)=W^*\alpha(WaW^*)W.$$
 Similarly one defines an endomorphic dilation of a semigroup
 $\{\Theta_t :t\in P\}$ of CP maps on $M$. If the semigroup is not
 discrete, one usually requires that certain continuity properties
 of the CP-semigroup would hold also for the endomorphism
 semigroup dilating it.

In \cite{Bh}, Bhat proved that every (unital) CP-semigroup $\{\Theta_t :
t\geq 0\}$ on the von Neumann algebra $B(H)$ can be dilated to a
(unital) semigroup of $^*$-endomorphisms on $B(K)$ for some larger
Hilbert space $K$. For general von Neumann algebras $M$ this was
proved by Bhat and Skeide in \cite{BS}. A different proof was
provided in \cite{MSQMP}. Both proofs used product systems of
correspondences but in a different way. In fact, the
correspondences in \cite{BS} are over $M$ while the
correspondences in \cite{MSQMP} are over $M'$. They are related by
 ``duality". (Since we shall not need it here, we will not
 elaborate on this concept of duality but  refer the reader to \cite{MSNP} or
 \cite{Sk}).

Since the methods of this paper will use some results and ideas
from \cite{MSQMP}, we shall now describe the approach taken there
(for a single CP map). Before we proceed, we note that, although
it was assumed in \cite{MSQMP} that the CP maps are unital, the
results we use here hold also for non unital maps.

Given a CP map $\Theta$ on a von Neumann algebra $M\subseteq
B(H)$, we write $M\otimes_{\Theta} H$ for the Hilbert space
obtained by the Hausdorff completion of the algebraic tensor
product $M\otimes H$ with respect to
$$\langle a\otimes h,b\otimes k\rangle= \langle h, a^*bk\rangle_H
\;,\;\;\; a,b\in M,\; h,k\in H.$$
A ``typical" element of $M\otimes_{\Theta}H$ will be written
$a\otimes_{\Theta}h$ and
there is a natural action of $M$ on this space where $a\in M$
sends $b\otimes_{\Theta}h$ to $ab\otimes_{\Theta}h$ (and we write
$a\otimes I_H$ for this operator). Now set

$$E_{\Theta}=\{X:H \rightarrow M\otimes_{\Theta}H :Xa=(a\otimes
I_H)X,\;a\in M \}.$$

As was shown in \cite[Proposition 2.5]{MSQMP},
this space is, in fact, a $W^*$-\\ correspondence over the von
Neumann algebra $M'$ (see Definition~\ref{corr}
) and there is a natural completely contractive representation associated to it.
The representation is $(\sigma,T_{\Theta})$ where $\sigma=id$, the
identity representation of $M'$, and $T(X)=W_{\Theta}^*X\in B(H)$ where
$W_{\Theta}:H \rightarrow M\otimes_{\Theta}H$ is defined by
$W_{\Theta}h=I\otimes_{\Theta}h$. One can check that $T_{\Theta}$
is an injective map (and so is $\sigma$).

To summarize, to every (contractive, normal) CP map on $M$ we
associated a pair $(E_{\Theta},(\sigma,T_{\Theta}))$ consisting of a
$W^*$-correspondence and a completely contractive representation
(and both $\sigma$ and
 $T_{\Theta}$ are injective).

This construction can be ``reversed".
Given a $W^*$-correspondence $E$ over $M'$ and a completely contractive
representation
$(\sigma,T)$ of $E$ on $H$ (such that the maps $\sigma$ and $T$ are injective),
 we can define a (contractive, normal)
CP map on $M$ by setting $\Theta_T(a)=\tilde{T}(I_E \otimes
a)\tilde{T}^*$, $a\in M$. (Here we use the Hilbert space
$E\otimes_{\sigma}H$ defined by the Hausdorff completion of the
algebraic tensor product with respect to $\langle \xi\otimes
h,\eta \otimes k\rangle=\langle h, \sigma(\langle
\xi,\eta\rangle)k\rangle$ and we let $\tilde{T}$ be the map
 $\tilde{T}:E\otimes_{\sigma}H \rightarrow H$ defined by
 $\tilde{T}(\xi \otimes_{\sigma}h)=T(\xi)h$ and $I_E \otimes a$ be
 the map sending $\xi \otimes_{\sigma}h$ to $\xi
 \otimes_{\sigma}ah$).

The two constructions are the inverse of each other up to
isomorphisms of pairs $(E,(\sigma,T))$ (that is, an isomorphism of
the correspondences that carries one representation to the other
one).
One direction of this statement is \cite[Corollary 2.23]{MSQMP}.
The other direction was proved in \cite{MSPS}.

Moreover, this bijection (between CP maps and pairs
$(E,(\sigma,T))$) carries $^*$-endomorphisms to representations
that are isometric (and vice versa). (See \cite[Proposition 2.21]
{MSQMP}).

The dilation of a single CP map can then be proved combining the
bijection described above with the dilation result for c.c.
representations (to isometric representations) in \cite[Theorem 3.3]{MS98}.
 For the details,
see \cite[Theorem 2.24]{MSQMP}.

In this paper we study to what extent we can apply these ideas to
product systems over $\mathbb{N}^2$ (in place of $\mathbb{N}$) and
a pair of commuting CP maps. The first result we need is the
dilation theorem for completely contractive representations of
product system over $\mathbb{N}^2$. This is achieved in
Theorem~\ref{dilation}. Applied to the case where $M=\mathbb{C}$
and each ``fiber" of the product system is $\mathbb{C}$, this
theorem yields Ando's Theorem  (for dilations of a pair of
commuting contractions to a pair of commuting isometries).
Since it is known that, in general, one cannot dilate
simultanuously a  commuting triple of contractions to a commuting triple
 of isometries (see \cite[Chapter
5]{Pau}), one cannot hope to have a general isometric dilation
result for representations of product systems over $\mathbb{N}^k$
for $k>2$.

A consequence of Theorem~\ref{dilation}
(Corollary~\ref{diltuples})
 is that two row contractions that, in some general sense, commute with
 each other, can be simultanuously dilated to two isometric row
 contractions preserving the commutation relation. (
Giving up the commutation relation, this result can be found in
\cite{Po2}.
 For a single
 row contraction, the dilation result was proved by Popescu in
 \cite{Po}).

Trying to extend the bijection described above (between CP maps
and pairs $(E,(\sigma,T))$) from the case $P=\mathbb{N}$ to the
case $P=\mathbb{N}^2$, one runs into a problem.
 It turns out that one
has to require that the two commuting CP maps $\Theta$ and $\Phi$  satisfy a
stronger condition (see Definition~\ref{stronglycommute}). A pair
of CP maps satisfying this condition is said to commute strongly.
the condition is needed so that we can find a product system
$X_{\Theta,\Phi}$ and a representation of it that will play the
role played by $(E_{\Theta},(\sigma,T_{\Theta}))$ in the case of a single
CP map $\Theta$ (see Proposition~\ref{rep}). Assuming that this
stronger condition holds, we establish the required bijection (see
Proposition~\ref{XTheta} and the discussion preceeding it).
 This bijection, together with Theorem~\ref{dilation}, implies
 that every pair of CP maps that commute strongly can be
 simultanuously dilated to a commuting pair of $^*$-endomorphisms.

However, the dilation result holds even if the CP maps commute but
not strongly. In order to prove it, we first have to show that
every pair of commuting CP maps can be ``realized" using some
representation of a product system over $\mathbb{N}^2$. This is
proved in
Proposition~\ref{noninj}. What we lose here (if the maps do not
commute strongly) is the uniqueness of the product system and the
representation. Proposition~\ref{noninj} is then applied to dilate
a general pair of commuting CP maps (Theorem~\ref{enddilation}).

As is shown in Proposition~\ref{ind}, knowing that the maps
commute strongly  has the additional advantage that, for each
 of the CP maps, the
correspondences associated with the map and with its dilation are
isomorphic. This was proved useful, for single CP maps, in
studying the index and the curvature of a CP map in \cite{MSCUR}.

\end{section}
\begin{section}{Preliminaries : Correspondences and
representations}
We begin by recalling the notions of a $C^*$-correspondence and a
\\
$W^*$-correspondence. For the general theory of Hilbert
$C^*$-modules which we use, we will follow \cite{L}. In
particular, a Hilbert $C^*$-module $E$ over a $C^*$-algebra $A$
will be a \emph{right} Hilbert $C^*$-module. We write
$\mathcal{L}(E)$ for the algebra of continuous, adjointable
$A$-module maps on $E$. It is known to be a $C^*$-algebra.

\begin{definition}\label{corr}
\begin{enumerate}
\item[(1)] A $C^*$-correspondence over a $C^*$-algebra $A$ is a
Hilbert $C^*$-module $E$ over $A$ endowed with the structure of a
left $A$-module via a $^*$-homomorphism $\varphi_E :A \rightarrow
\mathcal{L}(E)$.
\item[(2)] A Hilbert $W^*$-module over a von Neumann algebra $M$
is a Hilbert $C^*$-module over $M$ that is self dual (i.e., every
continuous $M$-module map from $E$ to $M$ is implemented by an
element of $E$).
\item[(3)] A $W^*$-correspondence over a von Neumann algebra $M$
is a Hilbert $W^*$-module $E$ that is a $C^*$-correspondence over
$M$ and the map $\varphi_E$ is a normal $^*$-homomorphism. (When
$E$ is a Hilbert $W^*$-module, $\mathcal{L}(E)$ is known to be a
von Neumann algebra \cite{Pa}).
\end{enumerate}
\end{definition}

When dealing with a specific $C^*$-correspondence $E$ it will be
convenient to write $\varphi$ (instead of $\varphi_E$) or even to
suppress it and write $a\xi$ or $a\cdot \xi$ for $\varphi(a)\xi$.

If $E$ and $F$ are $C^*$-correspondences over $A$, then the
balanced tensor product $E\otimes_A F$ is a $C^*$-correspondence
over $A$. It is defined as the Hausdorff completion of the
algebraic balanced tensor product with the internal inner product
given by
\be\label{tp}
\langle \xi_1 \otimes \eta_1,\xi_2\otimes \eta_2 \rangle=\langle
\eta_1, \varphi_F(\langle \xi_1,\xi_2 \rangle_E)\eta_2 \rangle_F
\ee
for all $\xi_1,\xi_2 \in E$ and $\eta_1,\eta_2 \in F$. The left
and right actions of $a\in M$ are defined by
\be\label{lrp}
\varphi_{E\otimes F}(a)(\xi \otimes \eta)b=\varphi_E(a)\xi \otimes
\eta b
\ee
for all $a,b\in M$, $\xi\in E$ and $\eta\in F$.

If $E$ and $F$ are $W^*$-correspondences over the von Neumann
algebra $M$, the tensor product $E\otimes_M F$ is understood to be
the self-dual extension (\cite{Pa}) of that Hausdorff completion.
The left and right actions are as in (\ref{lrp}) and, since
$\varphi_{E\otimes F}$ is now a normal $^*$-homomorphism,
$E\otimes_M F$ is a $W^*$-correspondence.

\begin{definition}\label{isomcor}
An isomorphism of $C^*$-correspondences (or
$W^*$-correspondences) $E$ and $F$  is a
surjective, bimodule map that preserves the inner products. We
write $E\cong F$ if such an isomorphism exists.
\end{definition}

If $E$ is a $C^*$-correspondence over $A$ and $\sigma$ is a
representation of $A$ on a Hilbert space $H$ (which is assumed to
be normal if $E$ is a $W^*$-correspondence) then
$E\otimes_{\sigma}H$ is the Hilbert space obtained as the
Hausdorff completion of the algebraic tensor product with respect
to $\langle \xi \otimes h, \eta \otimes k \rangle=\langle
h,\sigma(\langle \xi,\eta \rangle_E)k\rangle_H $. Given an
operator $X\in \mathcal{L}(E)$ and an operator $S\in \sigma(M)'$,
the map $\xi \otimes h \mapsto X\xi \otimes Sh $ defines a bounded
operator $X\otimes S$ on $E\otimes_{\sigma}H$. When $S=I_E$ and
$X=\varphi_E(a)$ (for $a\in A$) we get a representation of $A$ on
this Hilbert space. (If $E$ is a $W^*$-correspondence and $\sigma$
is a normal representation, so is $a\mapsto \varphi(a)\otimes
I_H$).
 We frequently
write $a\otimes I_H$ for $\varphi(a)\otimes I_H$.

\begin{definition}
\label{Definition1.12}Let $E$ be a $C^*$-correspondence over a $C^*$-algebra
 $A$. Then a completely contractive covariant
representation of $E$ (or, simply, a c.c. representation of $E$) on a
Hilbert space $H$ is a pair $(T,\sigma)$, where

\begin{enumerate}
\item[(1)] $\sigma$ is a $\ast$-representation of $A$ in $B(H)$.

\item[(2)] $T$ is a linear, completely contractive map from $E$ to
$B(H)$.


\item[(3)] $T$ is a bimodule map in the sense that $T(a\xi b)=\sigma(a)T(\xi
)\sigma(b)$, $\xi\in E$, and $a,b\in A$.
\end{enumerate}
If $A$ is a von Neumann algebra and $E$ is a $W^*$-correspondence,
we require also that
\begin{enumerate}
\item[(4)] $\sigma$ is a normal representation.
\end{enumerate}

\end{definition}

It should be noted that there is a natural way to view $E$ as an operator
space (by viewing it as a subspace of its linking algebra) and this defines
the operator space structure of $E$ to which the Definition
\ref{Definition1.12} refers when it is asserted that $T$ is completely contractive.

As we noted in the introduction and developed in \cite[Lemmas 3.4--3.6]{MS98}
and in \cite{MSNP}, if a completely contractive covariant representation,
$(T,\sigma)$, of $E$ in $B(H)$ is given, then it determines a contraction
$\tilde{T}:E\otimes_{\sigma}H\rightarrow H$ defined by the formula $\tilde
{T}(\eta\otimes h):=T(\eta)h$, $\eta\otimes h\in E\otimes_{\sigma}H$. The
operator $\tilde{T}$ satisfies
\begin{equation}
\tilde{T}(\varphi(\cdot)\otimes I)=\sigma(\cdot)\tilde{T}.\label{covariance}%
\end{equation}
In fact we have the following lemma from \cite[Lemma 2.16]{MSNP}.

\begin{lemma}
\label{CovRep}The map $(T,\sigma)\rightarrow\tilde{T}$ is a
bijection between all completely contractive covariant
representations $(T,\sigma)$ of $E$ on the Hilbert space $H$ and
contractive operators $\tilde{T}:E\otimes_{\sigma }H\rightarrow H$
that satisfy equation (\ref{covariance}). Given $\sigma$ and a
contraction $\tilde{T}$ satisfying the covariance condition
(\ref{covariance}), we get a  completely contractive covariant
representation $(T,\sigma)$ of $E$ on $H$ by setting $T(\xi)h:=\tilde{T}%
(\xi\otimes h)$.
\end{lemma}

\begin{remark}
\label{GenPowers}In addition to $\tilde{T}$ we also require the
\textquotedblleft generalized higher powers\textquotedblright\emph{\ }of
$\tilde{T}$. These are\emph{\ }maps$\;\tilde{T}_{n}:E^{\otimes n}\otimes
H\rightarrow H\;$defined by the equation$\;\tilde{T}_{n}(\xi_{1}\otimes
\ldots\otimes\xi_{n}\otimes h)=T(\xi_{1})\cdots T(\xi_{n})h$, $\xi_{1}%
\otimes\ldots\otimes\xi_{n}\otimes h\in E^{\otimes n}\otimes H$.
One checks easily that
 $\tilde{T}_n=\tilde{T}\circ
(I_E\otimes \tilde{T})\circ \cdots \circ (I_{E^{\otimes n-1}}
\otimes \tilde{T})$, $n>1$.

\end{remark}

\end{section}

\begin{section}{Representations of product systems and the
universal algebra}\label{ps} In the following we follow the
notation of Fowler (\cite{F}). Suppose $P$ is a left-cancellative,
countable, semigroup with an identity $e$ and $p:X \rightarrow P$
is a family of $C^*$-correspondences over $A$. Write $X_s$ for the
correspondence $p^{-1}(s)$ for $s\in P$ and $\varphi_s :A
\rightarrow \mathcal{L}(X_s)$ for the left action of $A$ on $X_s$.
We say that $X$ is a product system over $P$ if $X$ is a
semigroup, $p$ is a semigroup homomorphism and, for each $s,t \in
P \backslash \{e\}$, the map $(x,y) \in X_s \times X_t \rightarrow
X_{st}$ extends to an isomorphism $\theta_{s,t}$ of
correspondences from $X_s \otimes X_t$ onto $X_{st}$. We also
require that $X_e=A$ and that the multiplications $X_e\times X_s
\rightarrow X_s$ and $X_s\times X_e \rightarrow X_s$ are given by
the left and right actions of $A$ on $X_s$. The associativity of
the multiplication means that, for every $s,t,r \in P$,
\be\label{assoc} \theta_{st,r}(\theta_{s,t}\otimes
I_{X_r})=\theta_{s,tr}(I_{X_s}\otimes \theta_{t,r}). \ee

\begin{definition}\label{repn}
Suppose $H$ is a Hilbert space and $T:X\rightarrow B(H)$. Write
$T_s$ for the restriction of $T$ to $X_s$ and for $s=e$ write
$\sigma$ for $T_e$. We call $T$ (or $(\sigma,T)$) a completely
contractive representation of $X$ (and we write ``a c.c. representation") if
\begin{enumerate}
\item[(1)] For each $s$, $(\sigma,T_s)$ is a c.c. representation
of $X_s$ (as in Definition~\ref{Definition1.12}).
\item[(2)] $T(xy)=T(x)T(y)$ for all $x,y \in X$.
\end{enumerate}
Such a representation is said to be an isometric (or a Toeplitz)
representation if we also have
\begin{enumerate}
\item[(3)] $T(x)^*T(y)=\sigma(\langle x,y \rangle)$ whenever
$x,y\in X_s$ for some $s\in P$.
\end{enumerate}
\end{definition}

An important representation is the Fock representation. It is
defined as in \cite{F}. We write
$$\mathcal{F}(X) = \sum_{s \in P} \oplus X_s .$$
As mentioned in \cite{F}, this is a $C^*$-correspondence over $A$
with left action given by
$$\varphi_{\infty}(a)(\oplus x_s)=(\oplus \varphi_s(a)x_s).$$
We can define a representation $L$ of $X$ on $\mathcal{F}(X)$ by
setting
$$ L(x)(\oplus x_s) = \oplus (x\otimes x_s)\;,\;\; \oplus x_s \in
\mathcal{F}(X).$$

It is clear that $L$ is completely contractive. In fact, it is
completely isometric (i.e., $L_s$ is completely isometric for every $s\in P$).
 This can be seen even by considering the
restriction of $L(x)$ to $A\subseteq \mathcal{F}(X)$.

Note that, strictly speaking this is not what we defined as a
representation above (since $\mathcal{F}(X)$ is not a Hilbert space)
 but we can ``fix" it by representing
$\mathcal{L}(\mathcal{F}(X))$ on a Hilbert space.

As was shown in \cite[Proposition 2.8]{F}, the representation $L$ gives rise to a
$C^*$-representation of a certain $C^*$-algebra containing, for
 every $s\in P$, a copy
 of $X_s$ and the representation, restricted to this copy is equal
 to $L$. This $C^*$-algebra, $\mathcal{T}(X)$, (called \emph{the
 Toeplitz algebra} of $X$) has a universal property with respect
 to isometric (or Toeplitz) representations of $X$.

The next proposition shows that there is (a unique) operator
algebra $\mathcal{T}_+(X)$ which is universal with respect to c.c.
representations of $X$. The proof is standard and is omitted.

\begin{proposition}\label{universal}
Let $X$ be a product system over $P$ of $C^*$-correspondences over
$A$. Then there is a (closed) operator algebra $\mathcal{T}_+(X)$, called
the \emph{tensor algebra} of $X$, and a c.c. representation
$(i_A,i_X)$ of $X$ into $\mathcal{T}_+(X)$ such that
\begin{enumerate}
\item[(a)] $\mathcal{T}_+(X)$ is generated by the image of
$(i_A,i_X)$.
\item[(b)] For every c.c. representation $(\sigma,T)$ of $X$ on $H$,
there is a completely contractive representation $T\times \sigma$
of $\mathcal{T}_+(X)$ into $B(H)$ such that $(T\times \sigma)\circ
i_A=\sigma$ and $(T\times \sigma)\circ (i_X)_s = T_s$ (for $s\in
P$).

We shall refer to the maps $(i_A,i_X)$ as the \emph{universal
maps}.
\end{enumerate}
The triple $(\mathcal{T}_+(X), i_A,i_X)$ is unique up to a canonical
 completely isometric
isomorphism and, for every $s\in P$, $(i_X)_s$ is a complete isometry.
\end{proposition}

\begin{remark}
For $P=\mathbb{N}$, it follows from \cite{MS98} that
$\mathcal{T}_+(X)$ is the tensor algebra defined there. Hence, in
this case, it can be realized as a sulalgebra of
$\mathcal{L}(\mathcal{F}(X))$.
\end{remark}

Note that for $A=\mathbb{C}$ and a product system $X$ with
one-dimensional fibers (i.e., $X_s=\mathbb{C}$) and multiplication
induced from the multiplication of $P$, the algebra
$\mathcal{T}_+(X)$ is the algebra $OA(P)$ of Blecher and Paulsen
(\cite{BP}).

\end{section}

\begin{section}{Product systems over $ \mathbb{N}^2$}

Now we consider the case $P=\mathbb{N}^2$ (where
$\mathbb{N}=\{0,1,\ldots \}$) and prove a dilation result which
can be viewed as the analogue of Ando's dilation theorem (for two
commuting contractions).

We start with setting some notation. For $(m,n)\in \mathbb{N}^2$
and a product system of correspondences $X$ on $\mathbb{N}^2$, it
will be convenient to write $X(m,n)$ (instead of $X_{(m,n)}$) for
the fiber at $(m,n)$. If we set $E=X(1,0)$ and $F=(0,1)$, then
$X(m,n)$ is isomorphic to $E^{\otimes m}\otimes F^{\otimes n}$.
For convenience, we shall write $E^m$ for $E^{\otimes m}$ (and
similarly for $F$) and write $X(m,n)=E^m \otimes F^n$. (In other
words, we shall take this isomorphism to be the identity.) In the
notation of the previous section, this implies that
$\theta_{(m,0)(0,n)}=id$ and, more generally,
$\theta_{(k,0)(m,n)}$ and $\theta_{(k,l)(0,n)}$ are identity maps
(for $k,l,m,n \in \mathbb{N}$). Now, $X(m,n)$ is also isomorphic
to $F^n \otimes E^m$. This isomorphism will be written $t_{m,n}$,
so that

$$t_{m,n} :E^m \otimes F^n \rightarrow F^n \otimes E^m .$$

In fact, $t_{m,n}=\theta_{(0,n)(m,0)}^{-1}$ and we write $t$ for
$t_{1,1}$. Then, the associativity requirement enables one to
write each $t_{m,n}$ in terms of $t$. Straightforward computation
shows that we have
\be\label{t1n} t_{1,n}=(I_{F^{n-1}}\otimes
t)(I_{F^{n-2}}\otimes t \otimes I_F) \cdots (t\otimes I_{F^{n-1}})
\ee and \be\label{tmn} t_{m,n}=(t_{1,n}\otimes I_{E^{m-1}})(I_E
\otimes t_{1,n}\otimes I_{E^{m-2}})\cdots (I_{E^{m-1}}\otimes
t_{1,n}). \ee

Also, given an isomorphism $t:E\otimes F\rightarrow F\otimes E$,
we can define $t_{m,n}$ (using (\ref{t1n}) and (\ref{tmn})) and
use it to define $\theta_{(m,n)(k,l)}$ (for all $k,l,m,n \in
\mathbb{N}$) such that (\ref{assoc}) holds. Thus, defining a
product system over $\mathbb{N}^2$ amounts  to defining a triple
$(E,F,t)$ where $E$ and $F$ are $C^*$-correspondences over the
same $C^*$-algebra and $t:E\otimes F\rightarrow F\otimes E$ is an
isomorphism of correspondences.

Every completely contractive representation of $X$ on $H$ is determined
by its restrictions to $A$, to $E$ and to $F$. Thus we write such
a representation as a triple $(\sigma,T,S)$ where $T$ and $S$ are
the restrictions to $E$ and $F$ respectively.
The image of $x=\xi_1 \otimes \xi_2 \otimes \cdots \otimes \xi_m \otimes \eta_1
\otimes \eta_2 \otimes \cdots \otimes \eta_n \in E^m\otimes F^n$
under the representation would then be $T(\xi_1) \cdots
T(\xi_m)S(\eta_1)\cdots S(\eta_n)$.



Using Lemma~\ref{CovRep} and Remark~\ref{GenPowers}, we can write
the last expression as $\tilde{T}_m (I_{E^m} \otimes
\tilde{S}_n)(x)$.
We have
\be\label{tilde}
\tilde{T}_m (I_{E^m} \otimes \tilde{S}_n) :
X(m,n)\otimes_{\sigma} H=E^m \otimes F^n \otimes H \rightarrow
H
\ee

and, using condition (2) of Definition~\ref{repn},
we get the following ``commutation" relation
\begin{equation}\label{commutemn}
\tilde{T}_m (I_{E^m} \otimes \tilde{S}_n)=\tilde{S}_n
(I_{F^n} \otimes \tilde{T}_m)\circ (t_{m,n} \otimes I_H).
\end{equation}
For $m=n=1$ we have
\be\label{commute}
\tilde{T}(I_E\otimes \tilde{S})=\tilde{S}(I_F\otimes
\tilde{T})\circ (t\otimes I_H).
\ee
In fact, a tedious computation, using (\ref{tmn}), (which we omit)
shows that (\ref{commute}) implies (\ref{commutemn}) for all
$n,m \in \mathbb{N}$. Reversing the arguments, one also verifies
the following lemma.

\begin{lemma}\label{ccrep}
If $(\sigma,T)$ and $(\sigma,S)$ are completely contractive
representations of $E$ and $F$ respectively that satisfy
(\ref{commute}), then (\ref{tilde}) defines a (completely
contractive) representation of $X$.
\end{lemma}

\begin{remark}\label{CW}
So far we dealt with a product system of $C^*$-correspondences
over a $C^*$-algebra $A$. In Section~\ref{CP} we shall be
interested in a product system of $W^*$-correspondences over
a von Neumann algebra
$M$. For such a product system, a c.c. representation $T$ is
assumed to have the property that $\sigma$ ($=T_e$) is a normal
representation of $M$. (Note that then, using \cite[Remark
2.6]{MSNP}, each $T_s$ will, automatically, be continuous with
respect to the $\sigma$-topology on $X_s$ and the $\sigma$-weak
topology on $B(H)$).

\end{remark}

Now we discuss isometric dilations of completely covariant
representations. In the following we fix the product system $X$
and we use the notation set above.

\begin{definition}\label{dil}
Let $(\sigma,T,S)$ be a completely contractive covariant
representation of $X$ on a Hilbert space $H$. An isometric
dilation of $(\sigma,T,S)$ is an isometric representation $(\rho,V,U)$
of $X$ on a Hilbert space $K$ containing $H$, such that
\begin{enumerate}
\item[(1)] $H$ reduces $\rho$ and
$\rho(a)|H=P_H\rho(a)|H=\sigma(a)$, for all $a\in M$,
\item[(2)] $K\ominus H$ is invariant under each $V(\xi)$ and each
$U(\eta)$ (for $\xi \in E$, $\eta \in F$); that is,
$P_HV(\xi)|K\ominus H=P_HU(\eta)|K\ominus H=0$, and
\item[(3)] for all $\xi \in E$ and $\eta \in F$,
$P_HV(\xi)|H=T(\xi)$ and $P_HU(\eta)|H=S(\eta)$.
\end{enumerate}
We shall say that such a dilation is minimal in case the smallest
subspace of $K$ containing $H$ and invariant under every $V(\xi)$,
$\xi \in E$, and every $U(\eta)$, $\eta \in F$, is all of $K$.
\end{definition}

Note that, if $(\rho,V,U)$ is an isometric dilation of
$(\sigma,T,S)$ as above, then, for every $\xi_1,\xi_2, \ldots \xi_n \in E$
and $\eta_1,\eta_2, \ldots \eta_m \in F$,
$$P_H V(\xi_1)\cdots V(\xi_n)U(\eta_1) \cdots U(\eta_m)|H=
 T(\xi_1)\cdots T(\xi_n)S(\eta_1) \cdots S(\eta_m).$$
Also, a similar statement holds for all ``mixed" products; e.g.
$$ P_HV(\xi_1)U(\eta_1)V(\xi_2)|H=T(\xi_1)S(\eta_1)T(\xi_2).$$

\begin{theorem}\label{dilation}
Let $(\sigma,T,S)$ be a c.c. representation of $X$ on $H$ as
above. Then there is a Hilbert space $K$ containing $H$ and a
minimal
isometric representation $(\rho,V,U)$ of $X$ on $K$ that dilates
$(\sigma,T,S)$.

If $\sigma$ is non degenerate and $E$ and $F$ are essential (where
$E$ is essential if the subspace spanned by $\varphi(A)E$ is dense
in $E$), then $\rho$ is nondegenerate.

If $X$ is a product system of $W^*$-correspondences and $\sigma$
is assumed to be normal then $\rho$ is also a normal
representation.
\end{theorem}
\begin{proof}
We write $H_0$ for the Hilbert space $H$ together with the
representation $\sigma$ on it (we refer to it as a Hilbert module
over $A$) and define a sequence $\{H_k\}$ of Hilbert modules over
$A$ inductively by $H_{k+1}=H_k \oplus H_k^{(\infty)}$ where $
H_k^{(\infty)}$ is the direct sum of infinitely many copies of
$H_k$ (as Hilbert modules over $A$). We write $\sigma_k$ for the
representation of $A$ on $H_k$ and we think of $H_k$ as contained
in $H_{k+1}$ where the inclusion map sends $h$ to $h \oplus
0^{\infty}$. Also, given a correspondence $Y$ over $A$, we get an
inclusion $Y\otimes_{\sigma_k}H_k \subseteq
Y\otimes_{\sigma_{k+1}} H_{k+1}$.

The space $K$ that we need is
$$K= \sum_{(m,n)} \oplus (X(m,n)\otimes_{\sigma_{\max\{m,n\}}}H_{\max\{m,n\}}) .$$
There is a natural representation of $A$ on $K$. We shall write
$\rho$ for it (and we shall also write $\rho$ for its restriction
to various $\rho(A)$-invariant subspaces of $K$).
Note that $K=\sum_{l=0}^{\infty} K(l)$ where we write
\be
K(l)=\sum_{\max\{m,n\}=l} E^m \otimes F^n \otimes_{\sigma_l}H_l .
\ee
The dilation will constructed in several steps.

We first define $V_2:E\otimes K \rightarrow K$ and $U_2:F\otimes
K \rightarrow K$ by their restrictions to $E\otimes E^m\otimes F^n \otimes
H_{\max\{n,m\}}$ and $F\otimes E^m\otimes F^n \otimes
H_{\max\{n,m\}}$ as follows.

For a fixed $n\geq 0$, we write $q_0$ for the projection of $H_n$
onto $H=H_0$ (which is contained in $H_n$) and, for $h_0\in H$ and
$e\in E$, we set $D_T(e)h_0=\Delta_T(e\otimes h_0) \in E\otimes H
\subseteq E\otimes H_n$ where $\Delta_T=(I_{E\otimes H}-
\tilde{T}^*\tilde{T})^{1/2}\in B(E\otimes H)$. We then let $V_0:
E\otimes H_n \rightarrow H_n \oplus (E\otimes H_n)$ be defined by
\be
V_0(e\otimes h)=(T(e)q_0h)\oplus (D_T(e)q_0h \oplus (e\otimes
(I-q_0)h)). \ee Now, for $n=m=0$, we define $V_2|E\otimes H$ to be
$V_0$ (with $n=0$) and, for $m=0$ and $n>0$, we set
 $$V_2|E\otimes F^n \otimes_{\sigma_n}H_n=(I_{F^n \otimes
H_n}\oplus (t_{1,n}^{-1} \otimes I_{H_n}))\circ (I_{F^n}\otimes
V_0)\circ (t_{1,n}\otimes I_{H_n}).$$ Thus, for $n\geq 0$, $V_2$
maps $E\otimes F^n \otimes_{\sigma_n}H_n$ into
$(F^n\otimes_{\sigma_n}H_n)\oplus (E\otimes
F^n\otimes_{\sigma_n}H_n)$.

 Since $\norm{T(e)h_0 \oplus D_T(e)h_0}=\norm{\tilde{T}(e\otimes
 h_0) \oplus \Delta_T(e \otimes h_0)}=\norm{e\otimes h_0}$
, for $h_0\in H$, the map $V_0$ is an isometry. It is also
straightforward to check that $V_0$ is an $A$-module map (where
$a\in A$ acts on $H_n$ by $\sigma_n(a)$ and on $E\otimes H_n$ by
$\varphi_E(a)\otimes I_{H_n}$). Thus $I_{F^n}\otimes V_0$ is an
isometry from $F^n \otimes E\otimes H_n$ into $(F^n \otimes
H_n)\oplus (F^n \otimes E\otimes H_n)$. It follows that
 $V_2|E\otimes F^n \otimes_{\sigma_n}H_n$ is a composition of
 three isometries. Thus it is an isometry into $(F^n\otimes
 H_n)\oplus (E\otimes F^n\otimes H_n)$.

For $m>0$ we let $V_2|E\otimes E^m \otimes F^n \otimes
H_{\max\{n,m\}}$ be the inclusion map into $E^{m+1}\otimes F^n
\otimes H_{\max\{n,m+1\}}$ (where $E\otimes E^m$ is identified
with $E^{m+1}$ and $H_{\max\{n,m\}}$ is identified as a subspace
of $H_{\max\{n,m+1\}}$).

For different $n,m$ the ranges of $V_2|E\otimes E^m \otimes F^n \otimes
H_{\max\{n,m\}}$ are orthogonal to each other and, thus, it
follows that $V_2$ defines an isometry from $E\otimes K$ into $K$.

The definition of $U_2$ is similar. For $n=0$ we let
$$U_2|F\otimes E^m \otimes_{\sigma_m}H_m=(I_{E^m \otimes
H_m}\oplus (t_{1,m} \otimes I_{H_m}))\circ (I_{E^m}\otimes
U_0)\circ (t_{1,m}^{-1}\otimes I_{H_m})$$
where $U_0:F\otimes H_m \rightarrow H_m \oplus (F\otimes H_m)$ is
defined by
\be
U_0(f\otimes h)=(S(f)q_0h)\oplus (D_S(f)q_0h \oplus (f\otimes
(I-q_0)h)).
\ee
For $n>0$ we let $U_2|F\otimes E^m \otimes F^n \otimes
H_{\max\{n,m\}}$ be the map $t_{m,1}^{-1} \otimes I_{F^n} \otimes
I_{H_{\max\{n,m\}}}$ composed with the inclusion map of $E^m
\otimes F \otimes F^n \otimes H_{\max\{n,m\}}$
 into $E^{m}\otimes F^{n+1}
\otimes H_{\max\{n+1,m\}}$. Clearly, $U_2$ is an isometry from
$F\otimes K$ into $K$.

It is easy to check that we have, for $a\in A$,
\be\label{V}
V_2 (\varphi_E(a)\otimes I_{K})=\rho(a) V_2  \ee and  \be\label{U}
U_2 (\varphi_F(a)\otimes I_{K})=\rho(a) U_2. \ee
In general, the isometries $V_2,U_2$ do not necessarily satisfy a commutation
relation as in Equation (\ref{commute}). In fact, the maps
$V_2(I_E \otimes U_2)$ and $U_2(I_F \otimes V_2)(t \otimes
I_H)$ (where here $t$ is $t_{1,1}$), defined on
$E\otimes F \otimes H$, may differ. However, both maps map
$E\otimes F \otimes H$ into $H\oplus (E\otimes H)\oplus (F\otimes
H) \oplus (E\otimes F \otimes H)\subseteq K(0)\oplus K(1)$.
 For every $l\geq 0$, we write $P_l$ for the projection of
$K$ onto $K(l)$ (so that $P_0$ is the projection onto $H=H_0$).
 A simple computation shows that
$$P_0V_2(I_E \otimes U_2)=\tilde{T}(I_E \otimes
\tilde{S})=\tilde{S}(I_F\otimes \tilde{T})(t\otimes I_H)=
P_0U_2(I_F \otimes V_2)(t \otimes I_H).$$ Recall that $K(1)$ is
$(E\otimes H_1)\oplus (F\otimes H_1) \oplus (E\otimes F \otimes
H_1)$ and $P_1$ is the projection onto it. Write
 $G_1$ for the closure of the range of $P_1V_2(I_E \otimes U_2)$
 and $G_2$ for the closure of the range of $P_1U_2(I_F \otimes V_2)(t \otimes
I_H)$.

Since both maps $V_2(I_E \otimes U_2)$ and $U_2(I_F \otimes V_2)(t \otimes
I_H)$ are isometries, we have, for $\xi \in E\otimes F \otimes H$,
$\norm{P_1V_2(I_E \otimes U_2)\xi}^2=\norm{V_2(I_E \otimes
U_2)\xi}^2-\norm{P_0V_2(I_E \otimes
U_2)\xi}^2=\norm{\xi}^2-\norm{P_0U_2(I_F\otimes V_2)(t\otimes
I_H)\xi}^2=\norm{U_2(I_F\otimes V_2)(t\otimes
I_H)\xi}^2-\norm{P_0U_2(I_F\otimes V_2)(t\otimes
I_H)\xi}^2=\norm{P_1U_2(I_F\otimes V_2)(t\otimes
I_H)\xi}^2$.

Thus, the map sending $P_1V_2(I_E \otimes U_2)\xi$ to
$P_1U_2(I_F\otimes V_2)(t\otimes
I_H)\xi$ defines a unitary map from $G_1$ onto $G_2$.
We write $W(1)'$ for this map. For every $a\in A$, $\rho(a)$ maps
$G_i$ ($i=1,2$) into itself; so that we can view $G_i$ as a left
$A$-module. Moreover, the map $W(1)'$ is an isomorphism of
$A$-modules (i.e., the representations of $A$ associated with
$G_1$ and $G_2$ are equivalent).
Now write $\tau$ for the representation of $A$ on $(E\otimes H)
\oplus (F\otimes H) \oplus (E\otimes F \otimes H)$ (i.e., $\tau$
is the restriction of $\rho$ to this space) and
$\tau_{\infty}$ for the representation of $A$ on $(E\otimes H^{(\infty)})
\oplus (F\otimes H^{(\infty)}) \oplus (E\otimes F \otimes H^{(\infty)})$
(where $H^{(\infty)}=H_1 \ominus H$). Clearly,
$\tau_{\infty}$ is the sum of infinitely many copies of $\tau$.
Also write $\pi_i$ ($i=1,2$) for the representation of $A$ on $((E\otimes H)
\oplus (F\otimes H) \oplus (E\otimes F \otimes H))\ominus G_i$.
Then $\pi_i \leq \tau$ and, thus, $\pi_1 \oplus \tau_{\infty}
\cong \pi_2 \oplus \tau_{\infty}$. Write $W(1)'':((E\otimes H_1)
\oplus (F\otimes H_1) \oplus (E\otimes F \otimes H_1))\ominus G_1
\rightarrow ((E\otimes H_1)
\oplus (F\otimes H_1) \oplus (E\otimes F \otimes H_1))\ominus G_2$ for
the unitary implementing this equivalence and, forming $W(1):=W(1)'
\oplus
W(1)''$, we get a unitary operator on $K(1)$ that commutes with the
restriction of $\rho$ to $K(1)$. Also write $W(0)$ for the identity map on $H$.
 Then we have, for $\xi \in
E\otimes F \otimes H$,

$$W(1)V_2(I_E \otimes U_2)(I_E\otimes I_F \otimes W(0)^{-1})\xi=$$
$$U_2(I_F \otimes V_2)(t \otimes I_H)\xi.$$

Next, we shall define, inductively, unitary operators $W(k)$, on
$K(k)$, such that, writing $K[k]$ for the direct sum
$\sum_{l=0}^{k} \oplus K(l)$ and $W[k]$ for $\sum_{l=0}^k \oplus W(l)$
 for every $k\geq 0$, we get
\be\label{WV}
W[k+1]V_2(I_E \otimes U_2)(I_E\otimes I_F \otimes W[k]^{-1})\xi=
\ee
$$U_2(I_F \otimes V_2)(t \otimes I_{K[k]})\xi$$
for every $\xi \in E\otimes F \otimes K[k]$ . Once this is done, we write $W=\sum \oplus W(k)$,
 $\tilde{V}=WV_2$ and $\tilde{U}=U_2(I_F \otimes W^{-1})$ to get
\be\label{comm}
\tilde{V}(I_E \otimes \tilde{U})=
\tilde{U}(I_F \otimes
\tilde{V})(t \otimes I_K).
\ee

So we now assume that $W(l)$ has been defined for all $0\leq l
\leq k$. Write $G_1$ for the subspace $V_2(I_E \otimes
U_2)(I_E\otimes I_F \otimes W[k]^{-1})(E\otimes F \otimes K(k))$.
Since it is an isometric image of $E\otimes F \otimes K(k)$, it is
closed. Using the definition of $V_2$ and $U_2$ one can easily
check that $$G_1\subseteq \sum_{\max\{m,n\}=k+1} E^m \otimes F^n
\otimes H_k \subseteq K(k+1).$$ Similarly, write $G_2$ for the
(closed) subspace $U_2(I_F \otimes V_2)(t \otimes I_{K(k)})(E
\otimes F \otimes K(k))$ and note that $$G_2\subseteq
\sum_{\max\{m,n\}=k+1} E^m \otimes F^n \otimes H_k \subseteq
K(k+1).$$ The map $W(k+1)'$ sending $V_2(I_E \otimes
U_2)(I_E\otimes I_F \otimes W[k]^{-1})\xi$, in $G_1$ (for $\xi \in
E \otimes F \otimes K(k)$) to $U_2(I_F \otimes V_2)(t \otimes
I_{K(k)})\xi$, in $G_2$, is a unitary operator from $G_1$ onto
$G_2$ intertwining $\rho$. Now write $\pi_i$ for the restriction
of $\rho$ to $(\sum_{\max\{m,n\}=k+1} E^m \otimes F^n \otimes H_k)
\ominus G_i$, $\tau$ for the restriction of $\rho$ to
$\sum_{\max\{m,n\}=k+1} E^m \otimes F^n \otimes H_k$ and
$\tau_{\infty}$ for the restriction of $\rho$ to
$\sum_{\max\{m,n\}=k+1}\\ E^m \otimes F^n \otimes H_k^{\infty}$.
Now, argue as above (the case $k=0$) to find the unitary $W(k+1)$,
on $K(k+1)$, satisfying $$W(k+1)V_2(I_E \otimes U_2)(I_E\otimes
I_F \otimes W(k)^{-1})\xi= U_2(I_F \otimes V_2)(t \otimes
I_H)\xi=$$ $$ U_2(I_F \otimes V_2)(t \otimes I_{K(k)})\xi$$ for
each $\xi \in E\otimes F \otimes K(k)$. This, together with the
induction hypothesis,
 implies
(\ref{WV}) and, after setting $\tilde{V}=WV_2$ and
$\tilde{U}=U_2(I_F \otimes W)^{-1}$, we get (\ref{comm}).

Both $\tilde{V}$ and $\tilde{U}$ are isometries and it follows
from (\ref{V}) and (\ref{U}) and the fact that $W$ commutes with
$\rho(A)$, that, for $a\in A$,
$$ \tilde{V}(\varphi_E(a) \otimes I_K)=\rho(a)\tilde{V}$$ and
$$\tilde{U}(\varphi_F(a) \otimes I_K)=\rho(a)\tilde{U}.$$
Setting $V(\xi)k=\tilde{V}(\xi \otimes k)$ and
$U(\eta)k=\tilde{U}(\eta \otimes k)$ (for $\xi\in E$, $\eta \in F$
and $k \in K$), the triple
 $(\rho,V,U)$ defines an isometric
representation of $X$ on $K$. To see that it is a dilation of
$(\sigma,T,S)$ note that parts (1) and (2) of Definition~\ref{dil}
are easy to verify. To check part (3), fix $\xi \in E$ and $h\in
H\subseteq K$ and compute $P_HV(\xi)h=P_H\tilde{V}(\xi \otimes
h)=P_HWV_2(\xi \otimes h)=P_HWV_0(\xi \otimes
h)=P_H((T(\xi)q_0h)\oplus (D_T(\xi)q_0h \oplus (\xi \otimes
(I-q_0)h)))=T(\xi)h$. The computation for $U$ is similar.

The statement about the non degeneracy of $\rho$ is clear from its
definition. It is also clear that, if $E$ and $F$ are
$W^*$-correspondences over a von Neumann algebra $M$ and $\sigma$
is normal, so is $\rho$ (as both $\varphi_E$ and $\varphi_F$ are
assumed to be normal homomorphisms).

Finally, the dilation that we get in this way may not be minimal
but, restricting $(\rho,U,V)$ to the closed subspace of $K$
spanned by $H$ and by the vectors of the form $R_1R_2 \cdots
R_nh$, where $R_i \in U(F) \cup V(E)$ and $h\in H$, we get a
minimal isometric dilation.
\end{proof}

We now apply the theorem to obtain a dilation result for two
``commuting" row contractions. We note that, if one gives up the
commutativity condition in the next corollary, the dilation result
was obtained by Popescu in \cite{Po2}.

\begin{corollary}\label{diltuples}
Suppose $(T_1,\ldots T_n)$ and $(S_1,\ldots S_m)$ are an $n$-tuple
and an $m$-tuple of operators on a Hilbert space $H$ satisfying
\begin{enumerate}
\item[(i)] $\sum T_iT_i^*\leq I$ and $\sum S_jS_j^* \leq I$.
\item[(ii)] There is a unitary matrix $u=(u_{(i,j),(k,l)})$
(whose rows and columns are indexed by $\{1,2, \ldots n\}\times
\{1, \ldots m\}$) such that, for all $1\leq i \leq n$ and $1\leq j
\leq m$,
$$T_iS_j=\sum_{k,l}u_{(i,j),(k,l)} S_lT_k .$$
\end{enumerate}
Then there is a larger Hilbert space $K$, containing $H$, an
$n$-tuple of isometries $(V_1,\ldots V_n)$ in $B(K)$ and an
$m$-tuple of isometries $(U_1, \ldots U_m)$ in $B(K)$ such that
\begin{enumerate}
\item[(a)] $\sum V_iV_i^*\leq I$ and $\sum U_jU_j^* \leq I$ (that
is, in each tuple the isometries have pairwise orthogonal ranges).
\item[(b)] For every $1\leq i \leq n$ and $1\leq j
\leq m$,
$$V_iU_j=\sum_{k,l}u_{(i,j),(k,l)} U_lV_k .$$
\item[(c)] Each $V_i$ and $U_j$ leave $K\ominus H$ invariant.
\item[(d)] For every $1\leq i \leq n$ and $1\leq j
\leq m$, $P_HV_i|H=T_i$ and $P_HU_j|H=S_j$ (and, together with (c),
this implies that each product involving $V_i$'s and $U_j$'s dilates the
corresponding product with $T_i$'s and $S_j$'s).
\end{enumerate}
\end{corollary}
\begin{proof}
Let $M=\mathbb{C}$, $E=\mathbb{C}^n$, $F=\mathbb{C}^m$ (with
orthonormal bases $\{e_i\}$ and $\{f_j\}$ respectively)  and
$t:\mathbb{C}^n\otimes \mathbb{C}^m \rightarrow \mathbb{C}^m
\otimes \mathbb{C}^n$ be defined by $t(e_i\otimes
f_j)=\sum_{k,l}u_{(i,j)(k,l)} f_l \otimes e_k$. An $n$-tuple
$(T_1,\ldots T_n)$ satisfying $\sum T_iT_i^*\leq I$ defines a
completely contractive linear map $T:E\rightarrow B(H)$ by
$T(e_i)=T_i$. Similarly we define $S:F\rightarrow B(H)$ and (ii)
implies that they satisfy the commutation relation
(\ref{commute}).
 Letting $\sigma$ be the obvious
representation of $\mathbb{C}$ on $H$, we get a representation
$(\sigma,T,S)$ of the product system $X$ defined by $E$, $F$ and
$t$. Applying Theorem~\ref{dilation}, we get a Hilbert space $K$
and maps $V:E\rightarrow B(K)$ and $U:F\rightarrow B(K)$ defining
isometric representations (that dilate $T$ and $S$ respectively).
We now let $V_i$ be $V(e_i)$ and $U_j$ be $U(f_j)$.
The fact that these are isometric representations imply that
the operators $V_i$ and $U_j$ are all isometries. The rest of
(a)-(d) follows immediately.
\end{proof}

A special case of the following corollary (for $\alpha$ and
$\beta$ that are automorphisms) can be found in \cite{LM}.

\begin{corollary}\label{endom}
Let $\alpha$ and $\beta$ be commuting $^*$-endomorphisms of a
$C^*$-algebra $A$ that extend to the multiplier algebra $M(A)$ (as
commuting endomorphisms $\overline{\alpha}$ and
$\overline{\beta}$). Suppose $\sigma$ is a non degenerate
representation of $A$ on $H$ and $T_0$, $S_0$ are contractions in
$B(H)$ satisfying
\begin{enumerate}
\item[(i)] $\sigma(a)T_0=T_0\sigma(\alpha(a))$ and
$\sigma(a)S_0=S_0\sigma(\beta(a))$ for all $a\in A$, and
\item[(ii)] $T_0S_0=S_0T_0$.
\end{enumerate}
Then there is a Hilbert space $K$, containing $H$, a non
degenerate representation $\rho$ of $A$ on $K$ and partial
isometries $V_0$ and $U_0$ in $B(K)$ such that
\begin{enumerate}
\item[(1)] $\rho(a)V_0=V_0\rho(\alpha(a))$ and
$\rho(a)U_0=U_0\rho(\beta(a))$ for all $a\in A$,
\item[(2)] $V_0U_0=U_0V_0$,
\item[(3)] $U_0^*U_0=\overline{\rho}(\overline{\beta}(I))$ and
 $V_0^*V_0=\overline{\rho}(\overline{\alpha}(I))$,
 \item[(4)] $H$ reduces $\rho$ and $\rho(a)|H=\sigma(a)$,
 \item[(5)] $H$ is invariant for $U_0^*$ and for $V_0^*$, and
 \item[(6)]$P_HV_0|H=T_0$ and $P_HU_0|H=S_0$.
 \end{enumerate}
 \end{corollary}
 \begin{proof}
In the notation of Theorem~\ref{dilation}, let
$E=\overline{\alpha}(I)A=\overline{\alpha(A)A}$ and
$F=\overline{\beta}(I)A=\overline{\beta(A)A}$. The correspondence
structure of $E$ is defined by $\langle \xi,\eta
\rangle=\xi^*\eta$ and $\varphi_E(a)\xi b=\alpha(a)\xi b$, for
$a,b\in A$ and $\xi,\eta\in E$ (and similarly for $F$). Then one
can easily check that $E\otimes F$ is isomorphic to the
correspondence $\overline{\beta}\overline{\alpha}(I)A$ (via $\xi
\otimes \eta \mapsto \beta(\xi)\eta$) and $F\otimes E$ is
isomorphic to $\overline{\alpha}\overline{\beta}(I)A$. Combining
these isomorphisms, we get an isomorphism $t:E\otimes F
\rightarrow F\otimes E$ which can be written $t(\alpha(a_1)a_2
\otimes \overline{\beta}(I)b)=\overline{\beta}(I)\beta(a_1)\otimes
\overline{\alpha}\overline{\beta}(I)\beta(a_2)b$ for $a_1,a_2,b
\in A$.

A triple $(\sigma,T_0,S_0)$ satisfying (i) and (ii) defines a
representation $(\sigma,T,S)$ by setting
$T(\overline{\alpha}(I)a)=T_0\sigma(a)$ and
$S(\overline{\beta}(I)a)=S_0\sigma(a)$. Let $(\rho,V,U)$ be a
minimal isometric dilation. Then, for $a,b \in A$ and $h,k\in K$,
$\langle V(\overline{\alpha}(I)a)k,\\ V(\overline{\alpha}(I)b)g
\rangle= \langle \rho(b^*\overline{\alpha}(I)a)k,g \rangle =
\langle
\overline{\rho}(\overline{\alpha}(I))\rho(a)k,\rho(b)g\rangle$.
Thus, there is a partial isometry $V_0$ with
$V_0^*V_0=\overline{\rho}(\overline{\alpha}(I))$ satisfying
$V_0\rho(a)k=V(\overline{\alpha}(I)a)$. Similarly one defines
$U_0$ and properties (1) and (3) follow. Properties (4)-(6) follow
from the dilation properties and (2) follows from
Equation~\ref{commute} (for $\tilde{V}$ and $\tilde{U}$). We omit
the details.

 \end{proof}

\begin{remark}
For c.c. representations of a single $C^*$-correspondence
 it was shown in \cite[Theorem 4.4]{MS98} that commutant lifting
holds for the minimal isometric dilation. When $A=\mathbb{C}$,
this was proved in \cite[Theorem 3.2]{Po} generalizing the
commutant lifting theorem of Sz.-Nagy and Foias.
 It is known in the
classical case that the commutant lifting theorem of Sz.-Nagy and
Foias can be derived from Ando's dilation theorem. It is not hard
to see that Thoerem~\ref{dilation} can be used to give a different
proof of the commutant lifting theorem of \cite{MS98}. Since we
shall not use it in this paper, we omit the details.
\end{remark}

The following corollary shows that, when $P=\mathbb{N}^2$, the
universal tensor algebra $\mathcal{T}_+(X)$ of
Proposition~\ref{universal} is contained in the universal Toeplitz
$C^*$-algebra $\mathcal{T}(X)$ of \cite{F}.

\begin{corollary}\label{T+T}
Let $X$ be a product system of $C^*$-correspondences (over a
$C^*$-algebra $A$) with $P=\mathbb{N}^2$. Let $\mathcal{T}_+(X)$,
$i_A$ and $i_X$ be the universal tensor algebra and the universal
maps as in Proposition~\ref{universal}. Let $\mathcal{T}(X)$,
$k_A$ and $k_X$ be the universal Toeplitz algebra and the
universal maps as in \cite[Proposition 2.8]{F}. Then there is a
completely isometric homomorphism $$\Psi: \mathcal{T}_+(X)
\rightarrow \mathcal{T}(X)$$ such that $\Psi \circ i_A=k_A$ and
$\Psi \circ i_X=k_X$.
\end{corollary}
\begin{proof}
Write $\mathcal{B}$ for the norm-closed subalgebra of
$\mathcal{T}(X)$ generated by $k_A(A)$ and $k_X(X)$. We will show
that $(\mathcal{B},k_A,k_X)$ has the universal property (b) of
Proposition~\ref{universal}. Since it also satisfies (a), the
uniqueness of the universal algebra will complete the proof.

So suppose that $(\sigma,T)$ is a c.c. representation of $X$ on
$H$. It can be dilated to an isometric (i.e., Toeplitz)
representation $(\rho,V)$ on $K$. Then $(\rho,V)$ defines a
$C^*$-representation $V\times \rho$ of $\mathcal{T}(X)$ with
$(V\times \rho)\circ k_A=\rho$ and $(V\times \rho)\circ k_X=V$.
Set $\pi(b)=P_H(V\times \rho)(b)|H$ for $b\in \mathcal{B}$. Since
all $V(x)$ (for $x\in X$) and $\rho(a)$ (for $a\in A$) leave
$K\ominus H$ invariant, $P_H$ is a semiinvariant projection for
$(V\times \rho)(\mathcal{B})$ and, thus, $\pi$ is a completely
contractive representation of $\mathcal{B}$ on $H$. We also have,
for $a\in A$ and $x\in X$, $\pi(k_A(a))=P_H\rho(a)|H=\sigma(a)$
and $\pi(k_X(x))=P_HV(x)|H=T(x)$. Thus $\pi$ is $T\times \sigma$,
 completing the proof.

\end{proof}

\end{section}
\begin{section}{Commuting CP maps}\label{CP}
In this section we study commuting pairs of contractive, normal,
completely positive maps on von Nuemann algebras. The term ``CP
map" will always refer here to a contractive, normal, completely
positive map on a von Neumann algebra.

Let $\Theta$ and $\Phi$ be two normal CP maps on a given von
Neumann algebra $M$. We assume that $M\subseteq B(H)$ and consider
two Hilbert spaces defined as follows. On the algebraic tensor
product $M \otimes M \otimes H$ we define a sesquilinear form

$$\langle a_1\otimes b_1 \otimes h_1, a_2 \otimes b_2 \otimes h_2
\rangle = \langle h_1, \Theta (b_1^*
\Phi(a_1^*a_2)b_2)h_2\rangle.$$

We write $H_{\Phi,\Theta}$ (or
$M\otimes_{\Phi}M\otimes_{\Theta}H$)
 for the Hilbert space obtained by the Hausdorff completion of the
algebraic tensor product with respect to this semi inner product.
 A ``typical" element in
$H_{\Phi,\Theta}$ will be written $a\otimes_{\Phi}b
\otimes_{\Theta}h$. The Hilbert space $M \otimes_{\Phi} M
\otimes_{\Theta} H$ has a natural (normal) representation of $M$
on it. It is defined simply by $\lambda_{\Phi,\Theta}(c)(a\otimes
b \otimes h)=ca \otimes b \otimes h$, $c\in M$. It also has a
natural (normal) representation of $M'$ on it defined by
$\rho_{\phi,\Theta}(d)(a \otimes b \otimes h)=a \otimes b \otimes
dh$, $d\in M'$. We also write $c\otimes I_M \otimes I_H$ for
$\lambda_{\Phi,\Theta}(c)$ and $I_M \otimes I_M \otimes d$ for
$\rho_{\phi,\Theta}(d)$. We can, thus, think of $H_{\Phi,\Theta}$
as both a (normal) left $M$-module and a (normal) left
$M'$-module. Similarly, we define $H_{\Theta,\Phi}$ and view it as
a module over both $M$ and $M'$.

We now introduce a condition on the pair $(\Theta,\Phi)$ that is
stronger than the commutation relation $\Theta \Phi=\Phi \Theta$.
Its significance will be made clear later.

\begin{definition}\label{stronglycommute}
Given $\Theta$ and $\Phi$ as above, we say that they  commute
strongly
 if there is a unitary $u:H_{\Phi,\Theta} \rightarrow
H_{\Theta,\Phi}$ such that
\begin{enumerate}
\item[(i)] $u(a\otimes_{\Phi}I\otimes_{\Theta}h)=
a\otimes_{\Theta}I\otimes_{\Phi}h$ for $a\in M$ and $h\in H$.
\item[(ii)] $u(ca\otimes_{\Phi}b \otimes_{\Theta} h)=(c\otimes I_M \otimes
I_H)u(a\otimes_{\Phi}b \otimes_{\Theta} h)$ for $a,b,c\in M$ and
$h\in H$ (that is, $u$ intertwines the actions of $M$).
\item[(iii)] $u(a\otimes_{\Phi}b \otimes_{\Theta} dh)=(I_M\otimes
I_M \otimes d)u(a\otimes_{\Phi}b \otimes_{\Theta} h)$ for $a,b\in
M$, $d\in M'$ and $h\in H$ (that is, $u$ intertwines the actions
of $M'$).
\end{enumerate}
\end{definition}

\begin{remark}\label{commisom}
Note that, for $a\in M$ and $h \in H$, we have
$\norm{a\otimes_{\Phi}I\otimes_{\Theta}h}^2=\langle
h,\Theta(\Phi(a^*a))h \rangle$ while
$\norm{a\otimes_{\Theta}I\otimes_{\Phi}h}^2=\langle h,
\Phi(\Theta(a^*a))h\rangle$. Thus, the existence of a unitary $u$
satisfying (i)
 of Definition~\ref{stronglycommute} is equivalent to the assumption that $\Theta$ and $\Phi$
 commute. It follows that, if $\Theta$ and $\Phi$
commute strongly, then they commute. The converse is false, as we
shall see in Example~\ref{noncommstgly}.

\end{remark}

Given a von Neumann algebra $M\subseteq B(H)$ and a normal CP map
$\Theta: M\rightarrow M$, we write $$E_{\Theta}=\{X:H \rightarrow
M\otimes_{\Theta}H :\;Xa=(a\otimes I)X\;,\;a\in M\}.$$ (Recall
that $M\otimes_{\Theta}H$ was defined in the introduction). In
\cite{MSQMP} we wrote $\mathcal{L}_M(H,M\otimes_{\Theta}H)$ for it
and showed that it has a structure of a $W^*$-correspondence over
$M'$ (\cite[Proposition 2.5]{MSQMP}). In fact, the right action of
$d\in M'$ is given by $Xd=X\circ d$, the left action of $d$ is
$\varphi_{E_{\Theta}}(d)X=(I_M\otimes d)\circ X$ (where
$I_M\otimes d$ sends $a\otimes_{\Theta}h$ to
$a\otimes_{\Theta}dh$) and the inner product is $\langle
X_1,X_2\rangle=X_1^*X_2$, for $X_1,X_2\in E_{\Theta}$.
 We also defined (see \cite[Equation (2.7)
]{MSQMP}) the \emph{identity representation} of this
correspondence to be the pair $(\sigma,T_{\Theta})$ where $\sigma$
is the identity representation of $M'$ on $H$ and
$T_{\Theta}(X)=W_{\Theta}^*X$ (for $X\in E_{\Theta}$) where
$W_{\Theta}:H \rightarrow M\otimes_{\Theta}H$ is defined by
$W_{\Theta}h=I\otimes h$ (and, consequently,
$W_{\Theta}^*(a\otimes_{\Theta}h)=\Theta(a)h$). Note that
$T_{\Theta}:E_{\Theta} \rightarrow B(H)$ is an injective map
(because, for all $h,g \in H$ and $a\in M$, $\langle
W_{\Theta}^*Xa^*h,g \rangle= \langle Xa^*h,I\otimes_{\Theta}g
\rangle = \langle (I\otimes a^*)Xh,I\otimes_{\Theta}g \rangle
=\langle Xh,a\otimes_{\Theta}g \rangle$).

We also write (for normal CP maps $\Theta$ and $\Phi$)
$$E_{\Phi,\Theta}=\{Z:H \rightarrow H_{\Phi,\Theta}:
Za=(a\otimes_{\Phi}I \otimes_{\Theta}I)Z\;,\;a\in M\}.$$ Then
$E_{\Phi,\Theta}$ is a $W^*$-correspondence over $M'$ where the
right action is by composition, the left action is by
$\varphi(d)Z=(I\otimes I \otimes d)\circ Z$ and the inner product
is $\langle Z_1,Z_2 \rangle= Z_1^*Z_2$. Recall (\cite[Proposition
2.12]{MSQMP}) that the map $X\otimes Y \mapsto (I\otimes X)Y$ is
an isomorphism from the corresondence $E_{\Theta}
\otimes_{M'}E_{\Phi}$ onto the correspondence $E_{\Phi,\Theta}$.
We write $\Gamma_{\Phi,\Theta}$ for this isomorphism. Proposition
2.12 of \cite{MSQMP} also shows that there is an isometry $V$ from
 $E_{\Theta \Phi}$ into $E_{\Phi,\Theta}$  such that
$m:=V^*\Gamma_{\Phi,\Theta}$ is a coisometry mapping
$E_{\Theta}\otimes E_{\Phi}$ onto $E_{\Theta \Phi}$. Similarly,
one has a coisometry $n:E_{\Phi}\otimes E_{\Theta}\rightarrow
E_{\Phi \Theta}$.

\begin{remark}\label{comme}
It is easily seen from \cite[Proposition 2.12]{MSQMP} that (for
commuting maps $\Phi$ and $\Theta$) $\Theta$ and $\Phi$ strongly
commute if and only if the partial isometry $n^*m$ can be extended
to an isometry (of correspondences) from $E_{\Theta}\otimes
E_{\Phi}$ onto $E_{\Phi}\otimes E_{\Theta}$. In the case where
$M=B(H)$, these correspondences are Hilbert spaces (isomorphic to
Arveson's metric operator spaces, \cite{Arv2}) and the maps
commute strongly if and only if $\dim(Ker (m))=\dim(Ker (n))$.
\end{remark}

Using the remark above, the following lemma follows from
\cite[Proposition 2.14]{MSQMP}.

\begin{lemma}\label{stgcomm}
\begin{enumerate}
\item[(1)] If $\Theta$ and $\Phi$ are (normal) endomorphisms that
commute then they commute strongly.
\item[(2)] If $\Theta$ is a normal CP map and $\alpha$ is a normal
automorphism of $M$ that commutes with it then $\Theta$ and
$\alpha$
 commute strongly.
\item[(3)] If $\Theta$ is a normal CP map, $\alpha$ is a normal
automorphism of $M$ that commutes with $\Theta$ and $\Phi :=\Theta
\circ \alpha$ commutes with $\Theta$, then $\Theta$ and $\Phi$
commute strongly.
\end{enumerate}
\end{lemma}

\begin{example}\label{noncommstgly} There are pairs of
commuting normal CP maps that do not commute strongly.
\end{example}
Let $H$ be a Hilbert space, $P$ be a non trivial projection in
$B(H)$ and $S\in B(H)$ a coisometric map with $S^*S=P$ and such
that $S$ has some unit vector $k\in H$ with $S^*k=k$.
 Let
$\Theta: B(H) \rightarrow B(H)$ be the normal CP map
$\Theta(a)=\langle ak,k \rangle I_H$ and $\Phi:B(H)\rightarrow
B(H)$ be defined by $\Phi(a)=SaS^*$. Then, for $a\in B(H)$,
$\Phi(\Theta(a))=\Phi(\langle ak,k \rangle I)=\langle ak,k \rangle
SS^*=\langle aS^*k,S^*k \rangle I=\Theta(\Phi(a))$ so that the
maps commute and, in fact, $\Phi \circ \Theta=\Theta \circ
\Phi=\Theta$. A straightforward calculation shows that, for every
$a,b\in B(H)$ and $h\in H$, $ a \otimes_{\Phi} b \otimes_{\Theta}
h= aS^*bS \otimes_{\Phi} I \otimes_{\Theta} h$ in
$H_{\Phi,\Theta}$. Thus, $H_{\Phi,\Theta}$  is equal to the closed
subspace spanned by vectors of the form $c\otimes_{\Phi}
I\otimes_{\Theta} g$. On the other hand, if we choose $b\in B(H)$
and $h\in H$ such that $(I-P)bPh\neq 0$ and set
$x=I\otimes_{\Theta} (I-P)bP \otimes_{\Phi} h\in H_{\Theta,\Phi}$, then $x\neq 0$
and is orthogonal to the closed subspace of
$H_{\Theta,\Phi}$ spanned by the vectors of the form $c
\otimes_{\Theta} I \otimes_{\Phi} g$. This shows that the maps do
not commute strongly.

\vspace{5mm}

The importance of knowing whether two commuting normal CP maps
commute strongly follows from the next proposition. First, recall
that a (single) normal CP map on a von Neumann algebra $M$ always
``comes" from an (injective) representation of some
$W^*$-correspondence $E$. More precisely, given such CP map
$\Theta$ on $M\subseteq B(H)$, there is a $W^*$-correspondence $E$
over $M'$ and a completely contractive covariant representation
$(\sigma,T)$ of $E$ on $H$ (where $T$ is injective and
$\sigma=id$) such that

$$\Theta(a)=\tilde{T}(I_E \otimes a)\tilde{T}^*\;,\;\; a\in M.$$

(For the proof, see \cite[Corollary 2.23]{MSQMP}.) As the
following proposition shows, a similar statement holds for a
commuting pair of CP maps if and only if they commute strongly.

\begin{proposition}\label{rep}
Suppose $\Theta$ and $\Phi$ are commuting normal CP maps on
$M\subseteq B(H)$. Then the following are equivalent.
\begin{enumerate}
\item[(1)] $\Theta$ and $\Phi$ commute strongly.
\item[(2)] There is an isomorphism $t=t_{\Theta,\Phi}:E_{\Theta} \otimes_{M'}
E_{\Phi}\rightarrow E_{\Phi}\otimes_{M'}E_{\Theta}$ (defining a
product system $X_{\Theta,\Phi}$ over $\mathbb{N}^2$) such that
the identity representations $T_{\Theta}$ and $T_{\Phi}$ satisfy
\be\label{comm2}
 \tilde{T_{\Theta}} (I_{E_{\Theta}} \otimes
\tilde{T_{\Phi}})=\tilde{T_{\Phi}} (I_{E_{\Phi}} \otimes
\tilde{T_{\Theta}})\circ (t_{\Theta,\Phi} \otimes I_H) \ee
 (defining a
representation of the resulting product system such that, for
every $n,m$, $\tilde{T_{\Theta}}_m (I_{E^m} \otimes
\tilde{T_{\Phi}}_n)(I_{E^m} \otimes I_{F^n} \otimes a)(I_{E^m}
\otimes
\tilde{T_{\Phi}}_n)^*\tilde{T_{\Theta}}_m^*=\Theta^m(\Phi^n(a))
\;,\; a \in M$.)

\item[(3)] There is a product system $X(m,n)$
($(m,n)\in\mathbb{N}^2$) of $W^*$-correspondences over a von
Neumann algebra $N$ (with $E=X(1,0)$ and $F=X(0,1)$) and a
representation $(\sigma,T,S)$ of $X$ on $H$ such that $\sigma$ is
injective,
 $M=\sigma(N)'$, $T$ and $S$
are injective maps (of $E$ or $F$ into $B(H)$)
 and, for $a\in M$,
$\tilde{T}(I_E\otimes a)\tilde{T}^*=\Theta(a)$ and
$\tilde{S}(I_F\otimes a)\tilde{S}^*=\Phi(a)$.
\end{enumerate}
\end{proposition}
\begin{proof} We start by proving that (1) implies (2). Thus, we
assume that $\Theta$ and $\Phi$ commute strongly. It follows that
there is an isomorphism $u:H_{\Phi,\Theta}\rightarrow
H_{\Theta,\Phi}$ that maps $I\otimes_{\Phi}I \otimes_{\Theta} h$
to $I\otimes_{\Theta}I \otimes_{\Phi} h$ and satisfies the
conditions of Definition~\ref{stronglycommute}.
 Write $\Psi$ for the map
taking $Z\in E_{\Phi,\Theta}$ to $u \circ Z \in E_{\Theta,\Phi}$.
The fact that $u$ intertwines the representations of $M$ shows
that $\Psi(Z)$ is indeed in $E_{\Theta,\Phi}$. It is clearly an
isomorphism
 of $W^*$-modules. To see that it also intertwines the left
 actions
 of $M'$ on $E_{\Phi,\Theta}$ and on $ E_{\Theta,\Phi}$, we compute, for $d\in M'$,
 $\varphi(d)(u\circ Z)=(I_M\otimes I_M \otimes d)u\circ
 Z=u(I_M\otimes I_M\otimes d)Z=u \circ (\varphi(d)Z)$. Thus it is
 an isomorphism of $W^*$-correspondences.

  Recall that
$\Gamma_{\Phi,\Theta}$ is the isomorphism of
$E_{\Theta}\otimes_{M'}E_{\Phi}$ onto $E_{\Phi,\Theta}$ mapping
$X\otimes Y$ to $(I\otimes X)Y$ (\cite[Proposition 2.12]{MSQMP})
and write $t=t_{\Theta,\Phi}:E_{\Theta}\otimes_{M'}E_{\Phi}
\rightarrow E_{\Phi}\otimes_{M'}E_{\Theta}$
 for the isomorphism defined by
\be\label{t}
 t_{\Theta,\Phi}=\Gamma_{\Theta,\Phi}^{-1} \circ \Psi \circ
 \Gamma_{\Phi,\Theta}.
 \ee
We shall now turn to prove (\ref{commute}).

First, let $U_{\Theta}$ be the map from $M\otimes_{\Phi}H$ to
$M\otimes_{\Theta}M\otimes_{\Phi}H$ defined by
$U_{\Theta}(b\otimes_{\Phi}h)=I\otimes_{\Theta}b \otimes_{\Phi}h$.
Then $U_{\Theta}$ is a well defined contractive map and its
adjoint is $U_{\Theta}^*(a\otimes_{\Theta} b \otimes_{\Phi}
h)=\Theta(a)b\otimes_{\Phi}h$. Also, we have \be\label{wupsi}
W_{\Phi}^*U_{\Theta}^*u=W_{\Theta}^*U_{\Phi}^*. \ee To see this,
we compute, for  $h\in H$, $$
u^*U_{\Theta}W_{\Phi}h=u^*(I\otimes_{\Theta}I\otimes_{\Phi}h)=
I\otimes_{\Phi}I\otimes_{\Theta}h=U_{\Phi}W_{\Theta}h $$ and
(\ref{wupsi}) follows.

Note also that, for $X\in E_{\Theta}$,
$W_{\Phi}X^*(a\otimes_{\Theta}h)=I\otimes_{\Phi}(X^*(a\otimes_{\Theta}h))=
(I_M\otimes X^*)(I\otimes_{\Phi}a\otimes_{\Theta}h)=(I_M\otimes
X^*)U_{\Phi}(a\otimes_{\Theta}h)$. Thus, $U_{\Phi}^*(I_M\otimes
X)=XW_{\Phi}^*$ and, consequenly, for $X\in E_{\Theta}$ and $Y\in
E_{\Phi}$, \be\ U_{\Phi}^*\Gamma_{\Phi,\Theta}(X\otimes
Y)=U_{\Phi}^*(I_M\otimes X)Y= XW_{\Phi}^*Y. \ee It follows that,
for $h\in H$, \be\label{XY}
W_{\Theta}^*U_{\Phi}^*(\Gamma_{\Phi,\Theta}(X\otimes
Y))h=T_{\Theta}(X)T_{\Phi}(Y)h=\tilde{T_{\Theta}}(I \otimes
\tilde{T_{\Phi}})(X\otimes Y \otimes h). \ee Thus,
$\tilde{T_{\Theta}}(I\otimes \tilde{T_{\Phi}})(X\otimes Y \otimes
h)=
W_{\Theta}^*U_{\Phi}^*(\Gamma_{\Phi,\Theta}(X\otimes Y))h=
 W_{\Theta}^*U_{\Phi}^*\Psi^{-1}(\Gamma_{\Theta,\Phi}(t(X\otimes
 Y)))h=
 W_{\Phi}^*U_{\Theta}^*(\Gamma_{\Theta,\Phi}(
 t(X\otimes
 Y)))h=\tilde{T_{\Phi}}(I \otimes \tilde{T_{\Theta}})(t(X\otimes
 Y)\otimes h)=\tilde{T_{\Phi}}(I \otimes \tilde{T_{\Theta}})(t
 \otimes I_H))(X\otimes Y\otimes h)$. This proves (\ref{comm2}).

 Finally, the equation $\tilde{T_{\Theta}}_m (I_{E^m} \otimes
\tilde{T_{\Phi}}_n)(I_{E^m} \otimes I_{F^n} \otimes a)(I_{E^m}
\otimes
\tilde{T_{\Phi}}_n)^*\tilde{T_{\Theta}}_m^*=\Theta^m(\Phi^n(a))$
for $a\in M$ and arbitrary $m,n$ follows easily from the cases
$n=1, \;m=0$ and $m=1,\;n=0$. These, in turn, follow from
\cite[Corollary 2.23]{MSQMP}.

  This completes the proof
 that (1) implies (2). Since (2) obviously implies (3) (using
 \cite[Corollary 2.23]{MSQMP}), we now assume that (3) holds and
 turn to prove (1). As $\sigma$ is assumed to be injective and
 $M=\sigma(N)'$, we can replace $N$ by $\sigma(N)$ and assume
 $\sigma=id$ (and $N=M'$).

We start by defining the map $\Lambda_{\Theta,\Phi}:M
\otimes_{\Theta} M \otimes_{\Phi} H \rightarrow F \otimes_N E
\otimes_N H$ by
$$\Lambda_{\Theta,\Phi}(a\otimes_{\Theta}b\otimes_{\Phi}h)=(I_F\otimes
(I_E \otimes a)\tilde{T}^*b)\tilde{S}^*h$$ and the map
$\Lambda_{\Phi,\Theta}:M \otimes_{\Phi} M \otimes_{\Theta} H
\rightarrow E \otimes_N F \otimes_N H$ by
$$\Lambda_{\Phi,\Theta}(a\otimes_{\Phi}b\otimes_{\Theta}h)=(I_E\otimes
(I_F \otimes a)\tilde{S}^*b)\tilde{T}^*h.$$ We shall show that
these maps are (well defined, surjective) unitary maps and the map
\be \label{u} u:=\Lambda_{\Theta,\Phi}^{-1} \circ (t\otimes I_H)
\circ \Lambda_{\Phi,\Theta} \ee where $t:E\otimes_{M'}F\rightarrow
F\otimes_{M'}E$ is an isomorphism satisfying $\tilde{T}(I_E\otimes
\tilde{S})=\tilde{S}(I_F \otimes \tilde{T})(t\otimes I_H)$,
satisfies the conditions of Definition~\ref{stronglycommute}.

We first compute, for every $a,b,c,d\in M$ and $h,k \in H$,
$\langle
a\otimes_{\Theta}b\otimes_{\Phi}h,c\otimes_{\Theta}d\otimes_{\Phi}k
\rangle = \langle h, \Phi(b^*\Theta(a^*c)d)k \rangle = \langle h,
\tilde{S}(I_F \otimes b^*\tilde{T}(I_E \otimes
a^*c)\tilde{T}^*b)\tilde{S}^*k\rangle = \langle (I_F \otimes (I_E
\otimes a)\tilde{T}^*b)\tilde{S}^*h,(I_F \otimes (I_E \otimes
c)\tilde{T}^*d)\tilde{S}^*k \rangle $.

This shows that $\Lambda_{\Theta,\Phi}$ is a well defined
isometric map. Note also that, for $c \in M$, \be \label{c}
\Lambda_{\Theta,\Phi} \circ (c\otimes I_M \otimes I_H)=(I_F\otimes
I_E \otimes c)\Lambda_{\Theta,\Phi}. \ee To show that the map is
surjective, note first that the subspace of $F\otimes_N H$ spanned
by vectors of the form $(I_F \otimes b)\tilde{S}^*h$, for $h\in H$
and $b\in M$, is invariant under $I_F \otimes M$ and, thus, the
projection onto this subspace lies in $(I_F \otimes
M)'=\mathcal{L}(F) \otimes I_H$. Write $q\otimes I_H$ for it. If
$q\neq I$, there is some $\xi =(I-q)\xi \neq 0$ in $F$. But then,
for all $h,k \in H$ and $b\in M$, $0 = \langle \xi \otimes k, (I_F
\otimes b)\tilde{S}^*h \rangle=\langle \tilde{S}(\xi \otimes bk),h
\rangle =\langle S(\xi)bk,h \rangle$ contradicting the assumed
injectivity of $S$. Thus the closed subspace spanned by vectors of
the form $(I_F \otimes b)\tilde{S}^*h$ is all of $F\otimes H$.
Applying a similar argument to $T$ completes the proof that
$\Lambda_{\Theta,\Phi}$ is surjective.

Thus $\Lambda_{\Theta,\Phi}$ and $\Lambda_{\Phi,\Theta}$ are
unitary maps. Since $(t\otimes I_H)$ is unitary, so is $u$. Using
(\ref{c}), we find that $u$ intertwines the representations of
$M$. For $h\in H$, $u(I\otimes_{\Phi}I\otimes_{\Theta}h)=
\Lambda_{\Theta,\Phi}^{-1} \circ (t\otimes I_H) \circ
\Lambda_{\Phi,\Theta}(I \otimes_{\Phi}I\otimes_{\Theta}h)=
\Lambda_{\Theta,\Phi}^{-1} \circ (t\otimes I_H) (I_E \otimes
\tilde{S}^*)\tilde{T}^*h=\Lambda_{\Theta,\Phi}^{-1}(I_F \otimes
\tilde{T}^*)\tilde{S}^*h=I\otimes_{\Theta}I\otimes_{\phi}h$,
proving that $u$ satisfies  conditions (i) and (ii) of
Definition~\ref{stronglycommute}. To see that it satisfies
condition (iii), note that $\Lambda_{\Phi,\Theta}(I\otimes I
\otimes d)=(\varphi_E(d)\otimes I_F \otimes
I_H)\Lambda_{\Phi,\Theta}$
 (for $d\in M'$) and $t(\varphi_E(d)\otimes
 I_F)=(\varphi_F(d)\otimes I_E)t$.

\end{proof}

Now fix a von Neumann algebra $M\subseteq B(H)$.

Suppose $\Theta$ and $\Phi$ are normal CP maps on $M$ that
 commute strongly. Then, by
Proposition~\ref{rep}, we get a product system $X_{\Theta,\Phi}$
over $M'$, defined by $(E_{\Theta},E_{\Phi},t_{\Theta,\Phi})$,
and a representation $(id,T_{\Theta},T_{\Phi})$ of
$X_{\Theta,\Phi}$ on $H$ (and $T_{\Theta}$ and $T_{\Phi}$ are
injective maps). It will be convenient to refer to this
construction by $\tilde{X}$; that is,
$$\tilde{X}(\Theta,\Phi)=(X_{\Theta,\Phi},T_{\Theta},T_{\Phi}).$$

Conversely, suppose we start with a product system $X$ (of
$W^*$-\\correspondences over $M'$), defined by $(E,F,t)$ and suppose
 $(id,T,S)$ is a c.c. representation of $X$ on $H$
 (and $T$ and $S$ are injective maps). Then we get
 normal CP maps $\Theta$ and $\Phi$ on $M$ by setting
 $\Theta(a)=\tilde{T}(I_E\otimes a)\tilde{T}^*$ and
 $\Phi(a)=\tilde{S}(I_F\otimes a)\tilde{S}^*$ for $a\in M$.
  It follows from Proposition~\ref{rep} that $\Theta$ and $\Phi$
  commute strongly. We
 shall refer to this construction as $\tilde{\Theta}$; that is,
 $$\tilde{\Theta}(X,T,S)=(\Theta,\Phi).$$

Now, it follows from Proposition~\ref{rep} ((1) implies (2)) that
$$\tilde{\Theta}\circ \tilde{X}=id.$$
The following proposition shows that $\tilde{X}\circ
\tilde{\Theta}$ is an isomorphism. So that, up to isomorphisms of
product systems (more precisely, of product systems  with
representations), these two constructions are the inverses of each other.

\begin{proposition}\label{XTheta}
Let $M\subseteq B(H)$ be a von Neumann algebra.
Suppose $X$ is a product system  (of
$W^*$-correspondences over $M'$) defined by $(E,F,t)$ and suppose
 $(id,T,S)$ is a c.c. representation of $X$ on $H$
 (and $T$ and $S$ are injective maps). Let $\Theta$ and $\Phi$ be
 the normal CP maps defined on $M$ by
 $\Theta(a)=\tilde{T}(I_E\otimes a)\tilde{T}^*$ and
 $\Phi(a)=\tilde{S}(I_F\otimes a)\tilde{S}^*$ for $a\in M$.
Let $X_{\Theta,\Phi}$ be the product system constructed in the
proof of ``(1) implies (2)" of Proposition~\ref{rep} (so that it
is defined by $(E_{\Theta},E_{\Phi},t_{\Theta,\Phi})$ and
$t_{\Theta,\Phi}$ is as in (\ref{t})) and let
$(id,T_{\Theta},T_{\Phi})$ be the identity representation of
$X_{\Theta,\Phi}$.

Then there are (surjective) isomorphisms $w_E:E_{\Theta} \rightarrow E$
and $w_F:E_{\Phi} \rightarrow F$ such that
\begin{enumerate}
\item[(1)] $t\circ (w_E \otimes w_F)=(w_F\otimes w_E)\circ
t_{\Theta,\Phi}$, and
\item[(2)] $T\circ w_E=T_{\Theta}$ and $S\circ w_F= T_{\Phi}$.
\end{enumerate}
\end{proposition}

\begin{proof}
Let $v_E:M\otimes_{\Theta}H \rightarrow E\otimes H$ be defined by
$v_E(b\otimes h)=(I\otimes_{\Theta} b)\tilde{T}^*h$ and $v_F:M\otimes_{\Phi}H
\rightarrow F\otimes H$ is defined similarly (using $S$). The argument we gave
in the proof of Proposition`\ref{rep} to show that the map
$\Lambda_{\Theta,\Phi}$ is a unitary map shows also that $v_E$ and
$v_F$ are well defined unitary maps. (Note that this uses the injectivity
 of $T$ and $S$). It was shown in \cite[Theorem 2.14]{MSPS},
using the self duality of $E$,
that, for every $R\in E_{\Theta}$ one can find a (unique)
$w_E(R)\in E$ such that, for $\xi \in E$ and $h\in H$, $\langle
w_E(R),\xi \rangle h=R^*v_E^*(\xi \otimes h)$. It follows that,
for every $h\in H$ and $R\in E_{\Theta}$,
\be
w_E(R)\otimes h=v_ERh. \ee It is also shown there that $w_E$ is a
unitary, surjective, map from $E_{\Theta}$ onto $E$ and that part
(2) holds.

Now we turn to prove part (1). We first claim that, for every
$R\in E_{\Theta}$, $Y \in E_{\Phi}$ and $g\in H$, we have \be
\label{I} \Lambda_{\Phi,\Theta}(\Gamma_{\Phi,\Theta}(R\otimes
Y)g)=w_E(R)\otimes w_F(Y) \otimes g . \ee Recalling the definition
of $\Lambda_{\Phi,\Theta}$, we compute, for $a,b \in M$ and $h\in
H$, $\Lambda_{\Phi,\Theta}(a\otimes_{\Phi}b\otimes_{\Theta}h)=(I_E
\otimes I_F \otimes a)(I_E \otimes \tilde{S}^*)(I_E \otimes
b)\tilde{T}^*h=(I_E \otimes I_F \otimes a)(I_E \otimes
\tilde{S}^*)v_E(b\otimes_{\Theta} h)$. This holds for every
$b\otimes_{\Theta} h \in M\otimes_{\Theta}H$. In particular, it
holds with $Rf$ ($R\in E_{\Theta}$,$f\in H$) in place of
$b\otimes_{\Theta}h$. Thus $\Lambda_{\Phi,\Theta} (I_M\otimes R)(a
\otimes_{\Phi}f)= \Lambda_{\Phi,\Theta} (a\otimes_{\Phi}Rf)=(I_E
\otimes I_F \otimes a)(I_E \otimes \tilde{S}^*)v_E(Rf)=(I_E
\otimes I_F \otimes a)(I_E \otimes \tilde{S}^*)(w_E(R)\otimes f)=
w_E(R) \otimes (I_F\otimes a)\tilde{S}^*f=w_E(R)\otimes
v_E(a\otimes_{\Phi}f)$. Now write $Yg$ (for $Y\in E_{\Phi}$ and
$g\in H$) in place of $a\otimes_{\Phi}f$ to get
$\Lambda_{\Phi,\Theta} (I_M\otimes R)Yg=w_E(R)\otimes
v_E(Yg)=w_E(R)\otimes w_F(Y) \otimes g$. Since
$\Gamma_{\Phi,\Theta}(R\otimes Y)=(I_M \otimes R)Y$, this proves
(\ref{I}).

For $Z\in E_{\Phi,\Theta}$, $\Gamma_{\Phi,\Theta}^{-1}(Z)$ lies in
$E_{\Theta}\otimes E_{\Phi}$ and we can apply (\ref{I}) to it
(in place of $R\otimes Y$) and get
$\Lambda_{\Phi,\Theta}(Zg)=(w_E\otimes
w_F)(\Gamma_{\Phi,\Theta}^{-1}(Z))\otimes g$. Interchanging $\Phi$
and $\Theta$, we get $\Lambda_{\Theta,\Phi}(Gg)=(w_F\otimes
w_E)(\Gamma_{\Theta,\Phi}^{-1}(G))\otimes g$ for every $G\in
E_{\Theta,\Phi}$. Use it for $G=u \circ
\Gamma_{\Phi,\Theta}(R\otimes Y)$ (where
$u=\Lambda_{\Theta,\Phi}^{-1} \circ (t\otimes I_H)\circ
\Lambda_{\Phi,\Theta}$ as in (\ref{u})) to get
$\Lambda_{\Theta,\Phi}\circ u \circ \Gamma_{\Phi,\Theta}(R\otimes
Y)g=(w_F \otimes w_E)(t_{\Theta,\Phi}(R\otimes Y))\otimes g$. Now
use (\ref{u}) to get
\be \label{d}
(t\otimes I_H)\circ \Lambda_{\Phi,\Theta}\circ
(\Gamma_{\Phi,\Theta}(R\otimes Y))g=(w_F \otimes w_E)
(t_{\Theta,\Phi}(R\otimes Y))\otimes g.
\ee
But, using (\ref{I}), the left hand side of (\ref{d}) is equal to
$(t\otimes I_H)(w_E(R)\otimes w_F(Y) \otimes g)$.

This completes the proof of (1)

\end{proof}

\begin{proposition}\label{B(H)}
Let $\Theta$ and $\Phi$ be commuting normal CP maps on $B(H)$.
then they commute strongly if and only if
there are $n\leq \infty$ and $m\leq \infty$ and operators
$T_i,S_j$ in $B(H)$ ($1\leq i \leq n$, $1\leq j \leq m$) such
that
$$\Theta(a)=\sum_{i=1}^n T_i a T_i^*\;, \;a \in B(H), $$
$$\Phi(a)=\sum_{i=1}^m S_i a S_i^*\;, \;a \in B(H) $$
(where, if the sum is infinite, it is assumed to converge in the
weak operator topology) and $\{T_i\}$ and $\{S_j\}$ satisfy the
following conditions.
\begin{enumerate}
\item[(i)] $\sum T_iT_i^*\leq I$ and $\sum S_jS_j^*\leq I$.
\item[(ii)] ($l^2$-independence)  $\sum \alpha_i T_i\neq 0$
whenever $\alpha=\{\alpha_i\}\in l^2$ is nonzero (and similarly
for $\{S_j\}$).
\item[(iii)]  There is a unitary matrix $u=(u_{k,l}^{(i,j)})_{(i,j)(k,l)}$
 (whose rows and columns
are indexed by the set of pairs $(i,j)$ with $i\leq n$, $j\leq m$)
such that, for all $i,j$,
$$T_iS_j=\sum_{k,l} u_{k,l}^{(i,j)} S_l T_k. $$
\end{enumerate}
\end{proposition}
\begin{proof}
This is, in fact, a restatement of the equivalence of (1) and (3)
in Proposition~\ref{rep} for the case when $M=B(H)$.

\end{proof}

\begin{lemma}\label{pi}
Suppose $E$ and $F$ are $W^*$-correspondences over a von Neumann algebra
$N$ and $t:E\rightarrow F$ is a partial isometry in
$\mathcal{L}(E,F)$ that intertwines the left actions of $N$. (We shall refer to
such a map as a bimodule partial isometry). Then there  are projections $z_1$
and $z_2$
(in the center of $\mathcal{L}(E)\cap \varphi_E(N)'$ and the
center of
$\mathcal{L}(F)\cap \varphi_F(N)'$ respectively) and two
bimodule partial isometries $t_1,t_2$ in $\mathcal{L}(E,F)$ such
that
\begin{enumerate}
\item[(i)] $t_1^*t_1=z_1$ and $t_1t_1^*\leq z_2$ (so that we can view it as
a bimodule isometry from $z_1E$ into $z_2F$).
\item[(ii)] $t_2^*t_2\leq I_E - z_1$ and $t_2t_2^*=I_F-z_2$ (so
that we view it as a bimodule coisometry from $(I_E-z_1)E$ onto
$(I_F-z_2)F$).
\item[(iii)] $t_1$ extends $t_0z_1$ and $t_2$ extends
$t_0(I_E-z_1)$.
\item[(iv)] $(t_1+t_2)z_1=z_2(t_1+t_2)$.
\end{enumerate}
\end{lemma}
\begin{proof}
View $t_0$ as a partial isometry from $E\oplus F$ into $E\oplus F$
(by letting it be $0$ on $F$). Then it is a partial isometry in
the von Neumann algebra $R:=\mathcal{L}(E\oplus F) \cap
\varphi_{E\oplus F}(N)'$ (since it is a bimodule map).
 Apply the Comparison Theorem
(\cite[Theorem 6.2.7]{KR}) to the projections $f_1:=I_E-t_0^*t_0$
and $f_2:=I_F-t_0t_0^*$ to find a central projection $z$ in $R$
and partial isometries $v_1$ and $v_2$ in $R$ with
$v_1^*v_1=f_1z$, $v_1v_1^*\leq f_2z$, $v_2^*v_2\leq f_1(I-z)$ and
$v_2v_2^*=f_2(I-z)$. Finally, set $z_1=z I_E$, $z_2=z I_F$,
$t_1=t_0z + v_1$ and $t_2=t_0(I_E-z_1)+v_2$.

\end{proof}

\begin{lemma}\label{unitary}
Let  $E_0$ and $F_0$ be two $W^*$-correspondences over a von Neumann algebra
$N$ and
$t_0:E_0\otimes F_0 \rightarrow F_0\otimes E_0$ be a partial isometry
(of $W^*$-correspondences ; that is, it is an adjointable bimodule
map).
Then there is a partial isometry (of $W^*$-correspondences)
$t:E_0\otimes F_0 \rightarrow F_0\otimes E_0$ that extends $t_0$
and there
are $W^*$-correspondences $E$ and $F$ over $N$ containing $E_0$ and $F_0$
 respectively (as subcorrespondences), an isomorphism
(of correspondences)
$s:E\otimes F\rightarrow F\otimes E$ and projections $e_1$ and
$e_2$
such that ( writing $q_E$
 and $q_F$ for the projections of $E$ and $F$ onto $E_0$ and $F_0$
  respectively) we have
\begin{enumerate}
\item[(i)] $e_1$ lies in the center of $\mathcal{L}(E\otimes
F)\cap \varphi_{E\otimes F}(N)'$ and $e_2$ lies in the center of
 $\mathcal{L}(F\otimes
E)\cap \varphi_{F\otimes E}(N)'$.
\item[(ii)] $se_1(q_E\otimes q_F)=te_1(q_E \otimes q_F)=(q_F \otimes q_E)
te_1(q_E \otimes q_F)$ and this
map is an isometry from $e_1(E_0\otimes F_0)$ into $e_2(F_0\otimes
E_0)$.
\item[(iii)] $(I-e_2)(q_F \otimes q_E)s=t(I-e_1)(q_E\otimes q_F)=(q_F\otimes q_E)
t(I-e_1)(q_E\otimes q_F)$
and this map is a coisometry from $(I-e_1)(E_0 \otimes F_0)$ onto
$(I-e_2)(F_0\otimes E_0)$.
\item[(iv)] $te_1=e_2t$.

\end{enumerate}
\end{lemma}
\begin{proof}
Applying Lemma~\ref{pi} to $t_0$, we get projections $z_1$ (in the
center of $\mathcal{L}(E_0\otimes F_0)\cap \varphi_{E_0\otimes
F_0}(N)'$) and $z_2$ (in the
center of $\mathcal{L}(F_0\otimes E_0)\cap \varphi_{F_0\otimes
E_0}(N)'$) and partial isometries $t_1$ and $t_2$ (that are bimodule
 maps) satisfying
the conditions of that lemma. Write $t=t_1+t_2$. Then $t$ is a
partial isometry, $tz_1$ is an isometry from $z_1(E_0\otimes F_0)$
into $z_2(F_0\otimes E_0)$ and $t(I_{E_0\otimes F_0}-z_1)$ is a
coisometry from $(I_{E_0\otimes F_0}-z_1)(E_0\otimes F_0)$ onto
$(I_{F_0\otimes E_0}-z_2)(F_0\otimes E_0)$.

Let $E_1$ be a $W^*$-correspondence over $N$ that is isomorphic to
$E_0$ and let $F_1$ be isomorphic to $F_0$. These isomorphisms
induce (surjective) isomorphisms $\tau:E_0\otimes F_0\rightarrow
E_1 \otimes F_1$, $\theta:F_1\otimes E_1 \rightarrow F_0\otimes
E_0$
 and $\gamma:(E_0\otimes E_1)\oplus (F_1 \otimes
F_0) \rightarrow (E_1 \otimes E_0)\oplus (F_0 \otimes F_1)$. Write
$E=E_0\oplus F_1$ and $F= E_1 \oplus F_0$ and let $q_E$ and $q_F$
be the projections of $E$ onto $E_0$ and $F$ onto $F_0$
respectively. Also write $q_1$ for the projection $q_E\otimes q_F$
(from $E\otimes F$ onto $E_0 \otimes F_0$) and write $q_2$ for
$q_F \otimes q_E$. Clearly $q_1\in \mathcal{L}(E\otimes F)\cap
\varphi_{E\otimes F}(N)'$ and $q_2\in \mathcal{L}(F\otimes E)\cap
\varphi_{F\otimes E}(N)'$.

The isomorphism $s:E\otimes F \rightarrow F\otimes E$ will be
written matricially with respect to the decompositions $E\otimes
F=(E_0\otimes F_0)\oplus (F_1\otimes E_1)\oplus (E_0\otimes E_1
\oplus F_1 \otimes F_0)$ and $F\otimes E= (F_0\otimes E_0)\oplus
(E_1\otimes F_1)\oplus (E_1\otimes E_0 \oplus F_0 \otimes F_1)$ as
$$s=\left( \begin{array}{ccc} t & \theta-tt^*\theta & 0 \\ \tau -
\tau t^*t &
 \tau t^*\theta & 0 \\ 0 & 0 & \gamma
\end{array} \right) .$$
Clearly, $s$ is an isomorphism of correspondences.

Note that $\mathcal{L}(E_0\otimes F_0)\cap \varphi_{E_0\otimes
F_0}(N)'=q_1(\mathcal{L}(E\otimes F)\cap \varphi_{E\otimes
F}(N)'  )q_1$. Thus, there is a projection $e_1$ in the center of
$\mathcal{L}(E\otimes F)\cap \varphi_{E\otimes
F}(N)'$ such that $z_1=q_1e_1q_1$ (see \cite[Proposition 5.5.6 and
Corollary 5.5.7]{KR}). Similarly we get $e_2$ in the center of
 $\mathcal{L}(F\otimes E)\cap \varphi_{F\otimes
E}(N)'$ satisfying
$q_2e_2q_2=z_2$.



Since $t^*tz_1=t_1^*t_1z_1=z_1$, we see that $se_1q_1=sz_1=tz_1
=te_1q_1=q_2te_1q_1$ is an
isometry from $z_1(E_0\otimes F_0)$
into $z_2(F_0\otimes E_0)$ . This proves (ii) and a similar
argument works for (iii). Part (iv) here follows from part (iv) of
Lemma~\ref{pi}.

\end{proof}

\begin{proposition}\label{noninj}
Let $\Theta$ and $\Phi$ be two commuting normal CP maps on
$M\subseteq B(H)$. Then  there is a product system $X(m,n)$
($(m,n)\in\mathbb{N}^2$) of $W^*$-correspondences over the von
Neumann algebra $M'$
(with $E=X(1,0)$ and $F=X(0,1)$) and a representation
$(id,T,S)$ of $X$ on $H$ such that, for $a\in M$,
$\tilde{T}(I_E\otimes a)\tilde{T}^*=\Theta(a)$ and
$\tilde{S}(I_F\otimes a)\tilde{S}^*=\Phi(a)$.
\end{proposition}

\begin{proof}
We shall follow the idea of the proof of (1) implies (2) in
Proposition~\ref{rep} making changes when necessary. Since $\Theta$
and $\Phi$ commute, it follows from Remark~\ref{commisom} that
 there is a partial isometry $u_0$ in
$B(H_{\Phi,\Theta},H_{\Theta,\Phi})$ that is defined by the
formula
 $u_0(a\otimes_{\Phi}I
\otimes_{\Theta}h)=a\otimes_{\Theta}I\otimes_{\Phi}h$ (for $a\in
M$ and $h \in H$) and vanishes on the orthogonal complement of the
space spanned by the vectors $a\otimes_{\Phi}I \otimes_{\Theta}h$
, $a\in M$, $h \in H$. It is easy to check that $u_0$ intertwines
both the actions of $M$ on $H_{\Phi,\Theta}$ and on
$H_{\Theta,\Phi}$ and the actions of $M'$ on these spaces.
 We now write $\Psi_0$ for
the map taking $Z\in E_{\Phi,\Theta}$ to $u_0 \circ Z \in
E_{\Theta,\Phi}$. Then $\Psi_0$ is a partial isometry and a bimodule
map.
As in (\ref{wupsi}), we have
$W^*_{\Phi}U_{\Theta}^*=W_{\Theta}^*U_{\Phi}^*u_0^*$ and
$W_{\Theta}^*U_{\Phi}^*u_0^*u_0=W_{\Theta}^*U_{\Phi}^*$.
 Thus we have
 \be\label{Z}
W^*_{\Phi}U_{\Theta}^*Z=W_{\Theta}^*U_{\Phi}^*\Psi_0^*(Z)\;,\;
Z\in  E_{\Phi,\Theta}
\ee and
\be\label{R}
W^*_{\Theta}U_{\Phi}^*\Psi_0^*\Psi_0(R)=W^*_{\Theta}U_{\Phi}^*R\;,\;R\in
E_{\Theta,\Phi}.
\ee
 We now set
 $$t_0=\Gamma_{\Theta,\Phi}^{-1} \circ \Psi_0 \circ
 \Gamma_{\Phi,\Theta}$$
where $\Gamma_{\Phi,\Theta}$ is the isomorphism of
$E_{\Theta}\otimes_{M'}E_{\Phi}$ onto $E_{\Phi,\Theta}$ mapping
$X\otimes Y$ to $(I\otimes X)Y$. The map $t_0$ is a partial isometry
(of correspondences)
from $E_{\Theta}\otimes_{M'}E_{\Phi}$ to
$E_{\Phi}\otimes_{M'}E_{\Theta}$. Applying Lemma~\ref{unitary} to $t_0$, we
get $W^*$-correspondences $E$ and $F$ over $M'$, projections $e_1$ and $e_2$,
 a partial isometry $t$
extending $t_0$ and an isomorphism
$s:E\otimes_{M'}F \rightarrow F\otimes_{M'}E$ such that
$E_{\Theta}\subseteq E$, $E_{\Phi}\subseteq F$
and we have (writing $q_{\Theta}$ and $q_{\Phi}$ for the
projections onto $E_{\Theta}$ and $E_{\Phi}$ respectively),
\begin{enumerate}
\item[(i)] $e_1$ lies in the center of $\mathcal{L}(E\otimes
F)\cap \varphi_{E\otimes F}(M')'$ and $e_2$ lies in the center of
 $\mathcal{L}(F\otimes
E)\cap \varphi_{F\otimes E}(M')'$.
\item[(ii)] $se_1(q_{\Theta}\otimes q_{\Phi})=te_1(q_{\Theta} \otimes q_{\Phi})
=(q_{\Phi} \otimes q_{\Theta})te_1(q_{\Theta} \otimes q_{\Phi})$
 and this
map is an isometry from $e_1(E_{\Theta}\otimes E_{\Phi})$ into
$e_2(E_{\Phi}\otimes
E_{\Theta})$.
\item[(iii)] $(I-e_2)(q_{\Phi} \otimes q_{\Theta})s=t(I-e_1)
(q_{\Theta}\otimes q_{\Phi})=(q_{\Phi} \otimes q_{\Theta})t(I-e_1)
(q_{\Theta}\otimes q_{\Phi})$
and this map is a coisometry from $(I-e_1)(E_{\Theta} \otimes E_{\Phi})$ onto
$(I-e_2)(E_{\Phi}\otimes E_{\Theta})$.
\item[(iv)] $te_1=e_2t$.

\end{enumerate}

Define $\Psi=\Gamma_{\Theta,\Phi}\circ t\circ \Gamma_{\Phi,\Theta}^{-1}$.
Then $\Psi$ extends $\Psi_0$ (that is,
$\Psi_0=\Psi\Psi_0^*\Psi_0$)
and we have, using (\ref{Z}) and (\ref{R}),
$$W^*_{\Phi}U_{\Theta}^*Z=W_{\Theta}^*U_{\Phi}^*\Psi_0^*(Z)=
W_{\Theta}^*U_{\Phi}^*\Psi_0^*\Psi_0\Psi^*(Z)=
W_{\Theta}^*U_{\Phi}^*\Psi^*(Z)$$
for
$Z\in  E_{\Phi,\Theta}$.

 Now set $T=T_{\Theta}q_{\Theta}$,
$S=T_{\Phi}q_{\Phi}$ and $\sigma=id$.
 Then $(\sigma,T)$ and $(\sigma,S)$ are c.c. representations of $E$
 and $F$ respectively. (To see this, note
that $q_{\Theta}$ and $q_{\Phi}$ are bimodule maps since they
project onto subbimodules).

Note that (\ref{XY}) (in the proof of Proposition~\ref{rep}) still
holds whenever $X=q_{\Theta}(X)$ and $Y=q_{\Phi}(Y)$. So we can
  compute,
for $X\in E$, $Y\in F$ and $h\in H$,
$\tilde{S}(I\otimes \tilde{T})(e_2\otimes I_H)(Y\otimes X \otimes h)=
\tilde{T_{\Phi}}(I\otimes
\tilde{T_{\Theta}})(e_2\otimes I_H)
(q_{\Phi}(Y)\otimes q_{\Theta} (X) \otimes
h)=
W_{\Phi}^*U_{\Theta}^*(\Gamma_{\Theta,\Phi}e_2(q_{\Phi}(Y)\otimes
 q_{\Theta} (X)))h=W_{\Theta}^*U_{\Phi}^*\Psi^*(\Gamma_{\Theta,\Phi}e_2
 (q_{\Phi}(Y)\otimes
 q_{\Theta} (X)))h=
 W_{\Theta}^*U_{\Phi}^*(\Gamma_{\Phi,\Theta}(t^*(e_2(q_{\Phi}(Y)\otimes
q_{\Theta} (X))))h=
\tilde{T_{\Theta}}(I \otimes \tilde{T_{\Phi}})
 (t^*e_2(q_{\Phi}(Y)\otimes
q_{\Theta} (X))\otimes h)=\tilde{T_{\Theta}}(I \otimes
\tilde{T_{\Phi}})(t^*
 \otimes I_H))(e_2\otimes I_H)(q_{\Phi}(Y)\otimes q_{\Theta}(X)\otimes h)=
 \tilde{T_{\Theta}}(I\otimes \tilde{T_{\Phi}})(q_{\Theta}\otimes
 q_{\Phi}\otimes I_H)
 (t^* \otimes I_H)(e_2\otimes I_H) (q_{\Phi}(Y)\otimes q_{\Theta}(X)\otimes
 h)$.

 Using (ii) above, we find that $te_1(q_{\Theta} \otimes
 q_{\Phi})=e_2t(q_{\Theta} \otimes q_{\Phi})$ and, therefore,
 $e_1(q_{\Theta}\otimes q_{\Phi})t^*=(q_{\Theta}\otimes
 q_{\Phi})t^*e_2$

Thus $\tilde{S}(I\otimes \tilde{T})(e_2\otimes I_H)(Y\otimes X \otimes h)=
\tilde{T_{\Theta}}(I\otimes \tilde{T_{\Phi}})(e_1\otimes I_H)(q_{\Theta}\otimes
 q_{\Phi}\otimes I_H)
 (t^* \otimes I_H) (q_{\Phi}(Y)\otimes q_{\Theta}(X)\otimes
 h)=\tilde{T_{\Theta}}(I\otimes \tilde{T_{\Phi}})
 (e_1\otimes I_H)(q_{\Theta}\otimes
 q_{\Phi}\otimes I_H)(s^*\otimes I)(Y\otimes X \otimes h)=
 \tilde{T}(I \otimes \tilde{S})(e_1\otimes I_H)((s^*\otimes I)(Y\otimes X \otimes h)$.
 (Here we used (ii)  and the fact that
 $ \tilde{T}(I \otimes \tilde{S})= \tilde{T}(I \otimes
 \tilde{S})(q_{\Theta}\otimes q_{\Phi}\otimes I_H)$).

Since this holds for every $Y\otimes X \in F \otimes X$, it holds
for $s(X\otimes Y)$. Thus
\be \label{e1}
\tilde{T}(I \otimes \tilde{S})(e_1\otimes I_H)(X\otimes Y \otimes h)=
\tilde{S}(I\otimes \tilde{T})(e_2\otimes I_H)(s\otimes I_H)(Y\otimes X \otimes
h).
\ee
This dealt with the ``isometric" part. Now we turn to the
``coisometric" one and we compute, for $X,Y$ and $h$ as above,
$\tilde{T}(I\otimes \tilde{S})((I-e_1)\otimes I_H)(X\otimes Y \otimes h)=
\tilde{T_{\Theta}}(I\otimes
\tilde{T_{\Phi}})((I-e_1)\otimes I_H)(q_{\Theta}(X)\otimes q_{\Phi} (Y) \otimes
h)=
W_{\Theta}^*U_{\Phi}^*(\Gamma_{\Phi,\Theta}((I-e_1)(q_{\Theta}(X)\otimes
 q_{\Phi} (Y))))=
 W_{\Phi}^*U_{\Theta}^*\Psi(\Gamma_{\Phi,\Theta}((I-e_1)(q_{\Theta}(X)\otimes
 q_{\Phi} (Y))))=
 W_{\Phi}^*U_{\Theta}^*(\Gamma_{\Theta,\Phi}(t(I-e_1)(q_{\Theta}(X)\otimes
q_{\Phi} (Y))))h=
\tilde{T_{\Phi}}(I \otimes \tilde{T_{\Theta}})
 (t(I-e_1)(q_{\Theta}(X)\otimes
q_{\Phi} (Y))\otimes h)=\tilde{T_{\Phi}}(I \otimes
\tilde{T_{\Theta}})(t
 \otimes I_H))((I-e_1)\otimes I)(q_{\Theta}(X)\otimes q_{\Phi}(Y)\otimes h)=
 \tilde{S}(I \otimes \tilde{T})((I-e_2)\otimes I)((s\otimes I)(X\otimes Y \otimes h)$.
 Thus,
 \be\label{e2}
\tilde{T}(I\otimes \tilde{S})((I-e_1)\otimes I_H)(X\otimes Y \otimes h)=
\tilde{S}(I \otimes \tilde{T})((I-e_2)\otimes I_H)((s\otimes I_H)
(X\otimes Y \otimes h).
\ee
Adding up Equations (\ref{e1}) and (\ref{e2}), we get
\be\label{add}
\tilde{T}(I\otimes \tilde{S})=
\tilde{S}(I \otimes \tilde{T})(s\otimes I).
\ee
This shows that $(\sigma,T,S)$ is indeed a representation of the
system defined by $(E,F,s)$.

Finally, for $a\in M$,
$\tilde{T}(I_E\otimes
a)\tilde{T}^*=\tilde{T_{\Theta}}(q_{\Theta}\otimes I_H)(I_E
\otimes a)(q_{\Theta}\otimes I_H)\tilde{T_{\Theta}}^*=\tilde{T_{\Theta}}
(I_{E_{\Theta}}\otimes a)\tilde{T_{\Theta}}^*=\Theta(a)$ and a
similar computation applies to $\Phi$.

\end{proof}

\begin{definition}\label{dilcp}
Let $M\subseteq B(H)$ be a von Neumann algebra and let $\Phi$ and
$\Theta$ be two normal CP maps on $M$. An endomorphic dilation of
the pair $(\Phi,\Theta)$
 is a pair $(\alpha,\beta)$ of normal, commuting, $^*$-endomorphisms of a von
 Neumann algebra $R\subseteq B(K)$ and an
 isometry $W:H \rightarrow K$
such that
\begin{enumerate}
\item[(i)]  $M=W^*RW$,
\item[(ii)] $\alpha(WW^*) WW^*=\alpha(I)WW^*$ and $\beta(WW^*) WW^*=\beta(I)WW^*$,
\item[(iii)] $\Phi(a)=W^*\alpha(WaW^*)W$ and
$\Theta(a)=W^*\beta(WaW^*)W$ for all $a\in M$.
\end{enumerate}
\end{definition}

\begin{theorem}\label{enddilation}
Let $M\subseteq B(H)$ be a von Neumann algebra and $\Theta$ and
$\Phi$ be two commuting normal CP maps on $M$.
Then the pair $(\Theta,\Phi)$ has an endomorphic dilation.
\end{theorem}
\begin{proof}
Let $X$ and $(id,T,S)$ be as in Proposition~\ref{noninj}.
Using Theorem~\ref{dilation}, we find a Hilbert space $K$, an
isometric map $W:H \rightarrow K$ and an isometric representation
$(\rho,V,U)$ of $X$ on $K$ that dilates
$(id,T,S)$. Write $R=\rho(M')'\subseteq B(K)$ and let
$$ \alpha(b)=\tilde{V}(I_E\otimes b)\tilde{V}^*\;,\;\; b\in R $$ and
$$ \beta(b)=\tilde{U}(I_F\otimes b)\tilde{U}^* \;,\;\; b\in R.$$
Then, as was shown in Proposition~\ref{rep}, $\alpha$ and $\beta$
are two commuting, normal, CP maps on $R$. It follows from
\cite[Proposition 2.21]{MSQMP} that $\alpha$ and $\beta$ are
$^*$-endomorphisms of $R$.

Since $(\rho,V,U)$ is an isometric dilation of
$(id,T,S)$, we can write $W:H\rightarrow K$ for the isometric
embedding of $H$ into $K$ (so that $WW^*$ is the projection of $K$
onto $H$) and get
\be
W^*\tilde{V}(I_E \otimes W)=\tilde{T}
\ee and
\be
\tilde{V}^*WW^*=(I_E \otimes WW^*)\tilde{V}^*WW^*.
\ee
Thus, for $a\in M$,
\be
\Theta(a)=\tilde{T}(I_E\otimes a)\tilde{T}^*=
W^*\tilde{V}(I_E \otimes W)(I_E \otimes a)(I_E \otimes
W^*)\tilde{V}^*W=
\ee  \be =W^*\alpha(WaW^*)W
\ee and a similar equality holds for $\Phi$.
Also, $\alpha(I)WW^*=\tilde{V}\tilde{V}^*WW^*=
\tilde{V}(I_E \otimes WW^*)\tilde{V}^*WW^*=\alpha(WW^*)WW^*$. (And
similarly for $\beta$).

\end{proof}

\begin{remark}
 Consider the commuting CP maps of Example~\ref{noncommstgly}.
Let $S$ and $k$ be as assumed there. Note that $ \Theta(a)=\langle
ak,k\rangle I=\sum_{i=1}^{\infty} (e_i\otimes k^*)a(e_i \otimes
k^*)^* $ and
 $ \Phi(a)=SaS^*$
(where $\{e_i\}$ is a fixed orthonormal basis of $H$ and we write
$x\otimes y^*$ for the rank-one operator $(x\otimes y^*)h=\langle
h,y \rangle x$). Write $T_i=e_i \otimes k^*$. We can dilate
$\Theta$ to an endomorphism (of some $B(K)$) by using Popescu's
isometric dilation of the row contraction $\{T_i\}$ and we can
dilate $\Phi$ using an isometric dilation of the contraction $S$.
But these two endomorphisms will not commute. The proof of
Theorem~\ref{enddilation} shows that in order to get a commuting
pair, one needs to add zeroes to both families ($\{T_i\}$ and
$\{S\}$) before dilating them.
\end{remark}

In the case where $\Theta$ and $\Phi$ commute strongly we can say
more.

\begin{proposition}\label{ind}
If $\Theta$ and $\Phi$ are two normal CP maps on $M$ that commute
strongly then $(\Theta,\Phi)$ has an endomorphic dilation
$(\alpha,\beta)$ such that
\begin{enumerate}
\item[(1)] $E_{\Theta}$ is isomorphic (as a correspondence ) to
$E_{\alpha}$ and $E_{\Phi}$ is isomorphic to $E_{\beta}$.
\item[(2)] If $M$ is a semifinite factor and (using the
 terminology of \cite[Definition 4.9]{MSCUR}),
the index of $\Theta$ is finite, then so is the index of $\alpha$
and the two indices are equal. (Similar statement holds for $\Phi$
and $\beta$).
\end{enumerate}
\end{proposition}
\begin{proof}
Since $\Theta$ and $\Phi$ commute strongly, we can use, in the
proof of Theorem~\ref{enddilation} the product system
$X_{\Theta,\Phi}$ and the representation
$(id,T_{\theta},T_{\Phi})$ (instead of $X$ and $(id,T,S)$ that
were obtained from Proposition~\ref{noninj}). Writing $(\rho,V,U)$ for
an isomorphic dilation of this representation, we get $\alpha$ and $\beta$
 as in the proof of the theorem.

  In the notation of
the discussion preceeding Proposition~\ref{XTheta} we have
$(\alpha,\beta)=\tilde{\Theta}(X_{\Theta,\Phi},V,U)$. It follows
from Proposition~\ref{rep} (2) (since, as we have shown in
Lemma~\ref{stgcomm} (1), $\alpha$ and $\beta$ commute strongly)
that
$(\alpha,\beta)=\tilde{\Theta}(X_{\alpha,\beta},T_{\alpha},T_{\beta})$.
It now follows from Proposition~\ref{XTheta} that
$X_{\Theta,\Phi}$ is isomorphic to $X_{\alpha,\beta}$. In
particular, (1) follows. To be more precise, $E_{\Theta}$ and
$E_{\Phi}$ are correspondences over $M'$ while $E_{\alpha}$ and
$E_{\beta}$ are over $\rho(M')$. Thus the bimodule isomorphism and
the inner-product preservation are satisfied  ``up to $\rho$".
(Note that $\rho$ is injective). Part (2) follows immediately from
(1) and the fact that the index of a normal CP map $\Theta$
depends only on $E_{\Theta}$.

\end{proof}

\end{section}

\end{document}